\documentclass[12pt]{amsart}

\newtheorem{theorem}{Theorem}[section]
\newtheorem{lemma}[theorem]{Lemma}
\newtheorem{proposition}[theorem]{Proposition}

\newtheorem{corollary}[theorem]{Corollary}

\theoremstyle{definition}
\newtheorem{definition}[theorem]{Definition}
\newtheorem{notation}[theorem]{Notation}

\theoremstyle{remark}
\newtheorem{remark}[theorem]{Remark}

\def\myitemmarker{\P}
\newdimen\myindent
\setbox0=\hbox{\myitemmarker\ }
\myindent=\wd0

\def\myitem{\hangindent\myindent\hangafter-2\noindent
\llap{\hbox to\myindent{\myitemmarker\ }}}

\newcommand{\ab}{\allowbreak}
\newcommand{\alg}{\text{alg}}
\newcommand{\cA}{\mathcal{A}}

\newcommand{\CC}{{\mathbb C}}

\newcommand{\cP}{ {\mathcal P} }

\newcommand{\EE}{\mathrm{E}}
\newcommand{\ff}{\varphi}

\newcommand{\he}{{(\epsilon)}}

\newcommand{\kk}{{\mathrm k}}

\newcommand{\moeb}{\mathrm{M}\mathaccent"707F%
{\textrm{o}}\textrm{b}}

\newcommand{\NC}{NC}
\newcommand{\NN}{ {\mathbb N} }
\def\qed{{\unskip\nobreak\hfil\penalty50
 \hskip2em\hbox{}\nobreak\hfil 
\vbox to 7.7pt{\parindent0pt
\hsize7.7pt\hrule\vrule height7.3pt%
\hfill\vrule height7.3pt\hrule}
\parfillskip=0pt \finalhyphendemerits=0 \par}}

\newcommand{\RR}{ {\mathbb R} }
\newcommand{\SNC}{S_{NC}}
\newcommand{\SNCe}{S^{(\epsilon)}_{NC}}
\newcommand{\tr}{\mathrm{tr}}
\newcommand{\Tr}{\mathrm{Tr}}
\newcommand{\UU}{\mathcal{U}}
\newcommand{\Wg}{\mathrm{Wg}}
\newcommand{\ZZ}{ {\mathbb Z} }

\begin{document}
\allowdisplaybreaks
\title[Fluctuations of unitary random matrices]
{Second Order Freeness and Fluctuations of Random Matrices: \\
II. Unitary Random Matrices}

\author[J. A. Mingo]{James A. Mingo $^{(*)}$}
\address{Queen's University, Department of Mathematics and Statistics,
Jeffery Hall, Kingston, ON, K7L 3N6, Canada}
\email{mingo@mast.queensu.ca}
\thanks{$^*$ Research supported by Discovery Grants and a Leadership
Support Initiative Award from the Natural Sciences and Engineering
Research Council of Canada}

\author[P. \'Sniady]{Piotr \'Sniady $^{(\ddagger)}$}
\thanks{$\ddagger$ Research supported by State Committee for
Scientific Research (KBN) grant \mbox{2 P03A 007 23}, RTN
network: QP-Applications contract No.~HPRN-CT-2002-00279, and
KBN-DAAD project 36/2003/2004. The author is a holder of a
scholarship from the European Post-Doctoral Institute for
Mathematical Sciences.}
\address{Instytut Matematyczny, Uniwersytet Wroclawski, plac
Grun\-wald\-zki 2/4, 50-384 Wroclaw, Poland}
\email{Piotr.Sniady@math.uni.wroc.pl}

\author[R. Speicher]{\hbox{Roland Speicher $^{(*)(\dagger)}$}}
\thanks{$^\dagger$ Research supported by a Premier's
  Research Excellence Award from the Province of Ontario}
\address{Queen's University, Department of Mathematics and Statistics,
Jeffery Hall, Kingston, ON, K7L 3N6, Canada}
\email{speicher@mast.queensu.ca}

\begin{abstract}
We extend the relation between random matrices and free
probability theory from the level of expectations to the level
of fluctuations. We show how the concept of ``second order
freeness", which was introduced in Part I, allows one to
understand global fluctuations of Haar distributed unitary
random matrices. In particular, independence between the
unitary ensemble and another ensemble goes in the large $N$
limit over into asymptotic second order freeness. Two
important consequences of our general theory are: (i) we
obtain a natural generalization of a theorem of Diaconis and
Shahshahani to the case of several independent unitary
matrices; (ii) we can show that global fluctuations in
unitarily invariant multi-matrix models are not universal.
\end{abstract}
\date{}
\maketitle

\section{Introduction}

In Part I of this series \cite{MSp1} we introduced the concept of
second order freeness as the mathematical concept for dealing with
the large $N$ limit of fluctuations of $N\times N$-random matrices.
Whereas Voiculescu's freeness (of first order) provides the crucial
notion behind the leading order of expectations of traces, our
second order freeness is intended to describe in a similar way the
structure of leading orders of global fluctuations, i.e., of
variances of traces. In Part I we showed how fluctuations of
Gaussian and Wishart random matrices can be understood from this
perspective. Here we give the corresponding treatment for
fluctuations of unitary random matrices. Global fluctuations of
unitary random matrices have received much attention in the last
decade, see, e.g, the survey article of Diaconis \cite{D}.

Our main concern will be to understand the relation between unitary
random matrices and some other ensemble of random matrices which is
independent from the unitary ensemble. This includes in particular
the case that the second ensemble consists of constant (i.e.,
non-random) matrices. A basic result of Voiculescu tells us that on
the level of expectations, independence between the ensembles goes
over into asymptotic freeness. We will show that this result remains
true on the level of fluctuations: \textit{independence between
the ensembles implies that we have asymptotic second order
freeness between their fluctuations}.

\vskip2.5pt
Two important consequences of our investigations are the
following. 
\vskip2.5pt

\myitem
We get a generalization to the case of several
independent unitary random matrices of a classical result of
Diaconis and Shah\-shahani
\cite{DS}. Their one-dimensional case states that, for a unitary
random matrix $U$, the family of traces $\Tr(U^n)$ converge towards
a Gaussian family where the covariance between $\Tr(U^m)$ and
$\Tr(U^{*n})$ is given by $n\cdot \delta_{mn}$. In the case of
several independent unitary random matrices, one has to consider
traces in reduced words of these random matrices, and again these
converge to a Gaussian family, where the covariance between two such
reduced words is now given by the number of cyclic rotations which
match one word with the other. This result was also independently
derived by R\u{a}dulescu \cite{Rad} in the course of his
investigations around Connes's embedding problem.

\vskip2.5pt

\myitem
We can show that we do {\it not}\/ have
universality of fluctuations in multi-matrix models. For
unitarily invariant one-matrix models it was shown by
Johansson \cite{J} (compare also \cite{AJM}) that many random
matrix ensembles have the same fluctuations as the ensemble of
Gaussian random matrices. A main motivation for our
investigations was the expectation that many unitarily
invariant models of multi-matrix random ensembles should have
the same fluctuations as the ensemble of independent Gaussian
random matrices. However, our theory of second order freeness
shows that this is not the case.

The paper is organized as follows. In Section 2, we recall all the
necessary definitions and results around permutations, unitary
random matrices, and second order freeness. We will recall all the
relevant notions from Part I, so that our presentation will be
self-contained. However, for getting more background information on
the concept of second order freeness one should consult \cite{MSp1}.
In Section 3, we derive our main result about the asymptotic second
order freeness between unitary random matrices and another
independent random matrix ensemble. This yields as corollary that
independent unitary random matrices are asymptotically free of
second order, implying the above mentioned generalization of the
result of Diaconis and Shahshahani \cite{DS}. Section 4 shows how
our results imply the failure of universality of global
fluctuations in multi-matrix models.

\section{Preliminaries}

\subsection{The lattice of partitions}
For natural numbers $m,n\in\NN$ with $m<n$, we denote by $[m,n]$ the
interval of natural numbers between $m$ and $n$, i.e.,
$$[m,n]:=\{m,m+1,m+2,\dots,n-1,n\}$$ and $[m] = [1, m]$.
For a matrix $A=(a_{ij})_{i,j=1}^N$, we denote by $\Tr$ the
unnormalized and by $\tr$ the normalized trace,
$$\Tr(A):=\sum_{i=1}^N a_{ii},\qquad \tr(A):=\frac 1N \Tr(A).$$

We say that $A=\{A_1,\dots,A_k\}$ is a partition of the set $[1,n]$ if
the subsets $A_i$ are disjoint, non-empty,
and their union is equal to $[1,n]$. We call $A_1,\dots,A_k$ the
blocks of the partition $A$. For a permutation $\pi\in S_n$ we say that
a partition $A$ is $\pi$-invariant if $\pi$ leaves invariant each
block $A_i$. Let $1_n$ denote that partition of $[n]$ with one
block and $\cP(n)$ denote the partitions of $[n]$. Given
positive integers $m$ and $n$ let $1_{m,n}$ be the partition of
$[m+n]$ with the two blocks: $[m]$ and $[m+1, m+n]$.

If $A=\{A_1,\dots,A_k\}$ and $B=\{B_1,\dots,B_l\}$ are partitions of
the same set, we say that $A\leq B$ if for every block $A_i$ there
exists some block $B_j$ such that $A_i\subseteq B_j$. For a pair of
partitions $A,B$ we denote by $A\vee B$ the smallest partition $C$
such that $A\leq C$ and $B\leq C$. 

If we are considering classical random variables on a probability
space, then we denote by $\EE$ the expectation with respect to
the corresponding probability measure and by $\kk_r$ the
corresponding classical cumulants (as multi-linear functionals
in $r$ arguments); in particular,
$$\kk_1 (a )=\EE(a) \qquad\text{and}\qquad
\kk_2(a_1,a_2)=\EE( a_1a_2)-\EE(a_1) \EE(a_2).$$

If $a_1, \dots , a_n$ are random variables and $C = \{C_1, \dots , C_k
\}$ is in $\cP(n)$ we let 
$$
\EE_C(a_1, \dots , a_n) = \prod_{i=1}^k
\EE \Big( \mathop{\textstyle\prod}_{j \in C_i} a_j \Big)
$$
On the lattice $\cP(n)$  moments to cumulants are related by
the M\"obius function: $\moeb$. In particular the $n^{th}$
cumulant $\kk_n$ is given by
$$
\kk_n(a_1, \dots , a_r) = 
\sum_{C \in \cP(n)} \moeb(C, 1_r) \EE_C(a_1, \dots , a_n)
$$
where $\moeb(C, 1_n) = (-1)^{k-1} (k-1)!$ and where $k$ is the
number of blocks of $C$. We shall need the following formula for
the second cumulant of the product of random variables, see for
example \cite[3.2]{lehner}.

\begin{equation}\label{leonov}
\kk_2(a_1 \cdots a_m, b_1 \cdots b_n) =
\mathop{\sum_{\tau \in \cP(m+n)}}_{\tau \vee 1_{m,n} = 1_n}
\kk_\tau ( a_1, \dots a_m, b_1, \dots , b_n)
\end{equation}
The sum is over all partitions of $[m+n]$ which have at least
one block which connects the two blocks of $1_{m,n}$.

\subsection{Permutations}

We will denote the set of permutations on $n$ elements by $S_n$.  We
will quite often use the cycle notation for such permutations, i.e.,
$\pi=(i_1,i_2,\dots,i_r)$ is a cycle which sends $i_k$ to $i_{k+1}$
($k=1,\dots,r$), where $i_{r+1}=i_1$.

\subsubsection{Length function}

For a partition $\pi\in S_n$ we denote by $\#(\pi)$ the number of
cycles of $\pi$ and by $\vert\pi\vert$ the minimal number of
transpositions needed to write $\pi$ as a product of transpositions.
Note that one has
$$\vert\pi\vert+\#(\pi)=n\qquad\text{for all $\pi\in S_n$}.$$

\subsubsection{Non-crossing permutations}

Let us denote by $\gamma_n\in S_n$ the cycle
$$\gamma_n=(1, 2 ,\dots, n).$$
For all $\pi\in S_n$ one has that
$$n - 1 \leq \vert\pi\vert+\vert\gamma_n\pi^{-1}\vert.$$
If we have equality then we call $\pi$ \emph{non-crossing},
see \cite{Biane97} for the basic properties of non-crossing
partitions. Note that this is equivalent to
$$\#(\pi)+\#(\gamma_n\pi^{-1})=n+1.$$ If $\pi$ is
non-crossing, then so are $\gamma_n\pi^{-1}$ and $\pi^{-1}\gamma_n$;
the latter is called the \emph{(Kreweras) complement} of $\pi$.

We will denote the set of non-crossing permutations in $S_n$ by
$NC(n)$. Note that such a non-crossing permutation can be identified
with a non-crossing partition, by forgetting the order on the
cycles. There is exactly one cyclic order on the blocks of a
non-crossing partition which makes it into a non-crossing
permutation.

\subsubsection{Annular non-crossing permutations}

Fix $m,n\in\NN$ and denote by $\gamma_{m,n}$ the product of the two
cycles
$$\gamma_{m,n}=(1, 2,\dots, m)(m+1, m+2,\dots, m+n).$$
More generally, we shall denote by $\gamma_{m_1,\dots,m_k}$ the
product of the corresponding $k$ cycles.

We call a $\pi\in S_{m+n}$ \emph{connected} if the pair $\pi$ and
$\gamma_{m,n}$ generates a transitive subgroup in $S_{m+n}$. A
connected permutation $\pi\in S_{m+n}$ always satisfies
\begin{equation}\label{annular-NC}
m + n \leq \vert\pi\vert+\vert\gamma_{m,n}\pi^{-1}\vert.
\end{equation}
If $\pi$ is connected and if we have equality in that equation then
we call $\pi$ \emph{annular non-crossing}. Note that with $\pi$ also
$\gamma_{m,n}\pi^{-1}$ is annular non-crossing. Again, we call the
latter the \emph{complement} of $\pi$. Of course, all the above
notations depend on the pair $(m,n)$; if we want to emphasize this
dependency we will also speak about $(m,n)$-connected permutations
and $(m,n)$-annular non-crossing permutations.

We will denote the set of $(m,n)$-annular non-crossing permutations
by $\SNC(m,n)$. Again one can go over to annular non-crossing
partitions by forgetting the cyclic orders on cycles; however, in
the annular case, the relation between non-crossing permutation and
non-crossing partition is not one-to-one. Since we will not use the
language of annular partitions in the present paper, this is of no
relevance here.

Annular non-crossing permutations and partitions were introduced in
\cite{MN}; there, many different characterizations---in particular,
the one (\ref{annular-NC}) above in terms of the length
function---were given.

\subsection{A triangle inequality}

Let $\{A_1, \dots , A_k\}$ be a partition of $[n]$.
If, for $1\leq i\leq k$, $\pi_i$ is a permutation of the set $A_i$
we denote by $\pi_1\times \cdots\times \pi_k\in S_n$ the
concatenation of these permutations. We say that
$\pi=\pi_1\times\cdots\times \pi_k$ is a cycle decomposition if
additionally every factor $\pi_i$ is a cycle.

\begin{notation}
1) For $A\in\cP(n)$ we put $| A|:=n-\#(A)$ \\
2) For any $\pi \in S_n$ and any $\pi$-invariant $A\in\cP(n)$
we put
$$| (A,\pi) |:=2\vert A\vert-\vert\pi\vert.$$
\end{notation} 

\begin{lemma}\label{triangleinequality}
1 ) For all $A,B\in\cP(n)$ we have 
$$\vert A \vee B \vert\leq \vert A\vert+\vert B\vert$$
2) If $\pi$ and $\sigma$ are in $S_n$ and $A, B \in \cP(n)$ are
$\pi$ and $\sigma$ invariant respectively, then
$$\vert(A \vee B,\pi\sigma)\vert\leq \vert(A,\pi)\vert
+\vert(B,\sigma)\vert.$$
\end{lemma}

\begin{proof}
1) Each block of $B$ with $k$ points can glue together at
most $k$ blocks of $A$, thereby reducing the number of
blocks of $A$ by at most $k-1$. Thus $B$ can reduce by at
most $n-\#(B)$ the number of blocks of $A$. Hence the
difference between $\#(A)$ and $\#(A\vee B)$ cannot exceed $n-\#(B)$
and hence $$\#(A)-\#(A\vee B)\leq n-\#(B).$$ 
This is equivalent to our assertion.

2) We prove this, for fixed $\pi$ and $\sigma$ by induction
on $\vert A\vert+\vert B\vert$. The smallest possible
value of $\vert A\vert+\vert B\vert$ occurs when $\vert
A\vert=\vert\pi\vert$ and $\vert B\vert=\vert\sigma\vert$.
But then we have (since $A\vee B\geq \pi\sigma$)
$$2\vert A\vee B\vert-\vert\pi\sigma\vert
\leq \vert A\vee B\vert\leq \vert A\vert+\vert B\vert 
\eqno \textrm{by (1)}
$$
which is exactly our assertion for this case.
For the induction step note that we have just shown that
$$2\vert A\vee B\vert-\vert\pi\sigma\vert
\leq 2\vert A\vert - |\pi| + 2\vert B\vert - |\sigma| 
$$
when $|A| = |\pi|$ and $B| = |\sigma|$. Now one
only has to observe that if one increases $\vert A\vert$ (or
$\vert B\vert$) by 1 then $\vert A\vee B\vert$ can also
increase by at most 1. \end{proof}

\subsection{Haar distributed unitary random matrices and
            the Weingarten function}

In the following we will be interested in the asymptotics of special
matrix integrals over the group $\UU(N)$ of unitary $N\times
N$-matrices. We always equip the compact group $\UU(N)$ with its
Haar probability measure and accordingly distributed random matrices
we shall call \emph{Haar distributed unitary random matrices}. Thus
the expectation $\EE$ over this ensemble is given by integrating
with respect to the Haar measure.

The expectation of products of entries of Haar distributed unitary
random matrices can be described in terms of a special function on
the permutation group. Since such considerations go back to
Weingarten \cite{W}, Collins \cite{Col} calls this function the
Weingarten function and denotes it by $\Wg$. We will follow his
notation. In the following we just recall the relevant information
about this Weingarten function, for more details we refer to
\cite{Col,CS,X}.

We use the following definition of the Weingarten function. For
$\pi\in S_n$ and $N\geq n$ we put
$$\Wg(N,\pi)=
\EE \big(U_{11}\cdots
U_{nn}\overline{U}_{1\pi(1)}\cdots\overline{U}_{n\pi(n)}\big),
$$ where
$U=(U_{ij})_{i,j=1}^N$ is an $N\times N$ Haar distributed unitary
random matrix. Sometimes we will suppress the dependence on $N$ and
just write $\Wg(\pi)$. This $\Wg(N,\pi)$ depends on $\pi$ only
through its conjugacy class. General matrix integrals over the
unitary groups can be calculated as follows:
\begin{multline}
\EE\big(U_{i'_1j'_1}\cdots U_{i'_nj'_n}\overline U_{i_1j_1}
\cdots\overline U_{i_nj_n}\big) \\
= \sum_{\alpha,\beta\in S_n}
\delta_{i_1i'_{\alpha(1)}}\cdots \delta_{i_ni'_{\alpha(n)}}
\delta_{j_1j'_{\beta(1)}}\cdots\delta_{j_nj'_{\beta(n)}}
\Wg(\beta\alpha^{-1}).
\end{multline}

The Weingarten function is a quite complicated object, and its full
understanding is at the basis of questions around Itzykson-Zuber
integrals. For our purposes, only the behaviour of leading orders in
$N$ of $\Wg(N,\pi)$ is important. One knows (see, e.g.,
\cite{Col,CS}) that the leading order in $1/N$ is given by
$\vert\pi\vert+n$ and increases in steps of 2.

Let us use the following notation for the first two orders ($\pi\in
S(n)$):
$$\Wg(N,\pi)=\mu(\pi)N^{-(\vert \pi\vert +n)}+\phi(\pi)
N^{-(\vert\pi\vert+n+2)}+O\bigl(N^{-(\vert\pi\vert+n+4)}\bigr).$$

One knows that $\mu$ is multiplicative with respect to the cycle
decomposition, i.e.,
$$\mu(\pi_1\times\pi_2)=\mu(\pi_1)\cdot\mu(\pi_2).$$
The important part of the second order information is
contained in the leading order of
$\Wg(\pi_1\times\pi_2)-\Wg(\pi_1)\Wg(\pi_2)$, which is given by
$\mu_2(\pi_1,\pi_2) N^{-(|\pi_1| + |\pi_2| + m + n +2)}$ for
$\pi_1 \in S_m$ and $\pi_2 \in S_n$ and where
$$\mu_2(\pi_1,\pi_2):=
\phi (\pi_1\times\pi_2)-\mu(\pi_1)\phi(\pi_2)-
\phi(\pi_1)\mu(\pi_2).$$ Note that we have
$$\mu_2(\pi_1,\pi_2)=\mu_2(\pi_2,\pi_1).$$
Collins \cite{Col} has general counting formulas for the calculation
of $\mu$ and $\mu_2$ (and also higher order analogues); however, a
conceptual explanation of $\mu_2$ seems still to be missing. $\mu$
is the M\"obius function of the lattice of non-crossing
partitions (thus determined by Catalan numbers), and this fact
is quite well understood by the relation between $\mu$ and
asymptotic freeness of unitary random matrices. In a similar
way, one should get a conceptual understanding of $\mu_2$ by
the relation with second order freeness. In the present paper
we will not pursue further this direction, but we will come
back to it in forthcoming investigations. Here we will not
rely on the concrete values of
$\mu$ or $\mu_2$, but will only use the  basic
properties mentioned above.

\subsection{Second order freeness}

In \cite{MSp1}, we introduced the concept of second order freeness
which is intended to capture the structure of the fluctuation
functionals for random matrices arising in the limit $N\to\infty$,
in the same way as the usual freeness captures the structure of the
expectation of the trace in the limit. We recall the relevant
notations and definitions.

\begin{definition}
A \emph{second order non-commutative probability space}
$(\cA,\ab \ff_1,\ff_2)$ consists of a unital algebra $\cA$, a
tracial linear functional
$$\ff_1:\cA\to\CC \qquad\text{with}\qquad
\ff(1)=1$$ and a bilinear functional
$$\ff_2:\cA\times\cA\to\CC,$$
which is tracial in both arguments and which satisfies
$$\ff_2(a,1)=0=\ff_2(1,b)\qquad\text{for all $a,b\in\cA$.}$$
\end{definition}

\begin{notation}
Let unital subalgebras $\cA_1,\dots,\cA_r\subset\cA$ be given.
\\
1) We say that a tuple $(a_1,\dots,a_n)$ ($n\geq 1$) of elements
from $\cA$ is \emph{cyclically alternating} if, for each $k$, we
have an $i(k)\in\{1,\dots,r\}$ such that $a_k\in \cA_{i(k)}$ and, if
$n\geq 2$, we have $i(k)\not= i(k+1)$ for all $k=1,\dots,n$. We
count indices in a cyclic way modulo $n$, i.e., for $k=n$ the above
means $i(n)\not=i(1)$. Note that for $n=1$, we do not impose any
condition on neighbours.
\\
2) We say that a tuple $(a_1,\dots,a_n)$ of elements from $\cA$ is
\emph{centered} if we have
$$\ff_1(a_k)=0\qquad\text{for all $k=1,\dots,n$.}$$
\end{notation}

\begin{definition}\label{free}
Let $(\cA,\ff_1,\ff_2)$ be a second order non-commutative
probability space. We say that unital subalgebras
$\cA_1,\dots,\cA_r\subset \cA$ are \emph{free with respect to
$(\ff_1,\ff_2)$} or \emph{free of second order}, if they are free
(in the usual sense \cite{VDN}) with respect to $\ff_1$ and if the
following condition for $\ff_2$ is satisfied. Whenever we have, for
$n,m\geq 1$, tuples $(a_1,\dots,a_n)$ and $(b_m,\dots, b_1)$ from
$\cA$ such that both are centered and cyclically alternating then we
have
\begin{enumerate}
\item If $n\not=m$, then
$$\ff_2(a_1\cdots a_n, b_m\cdots b_1)=0.$$
\item If $n=m=1$ and $a\in\cA_i$, $b\in\cA_j$, with $i\not= j$,
then
$$\ff_2(a,b)=0.$$
\item If $n=m\geq2$, then
$$\ff_2(a_1\cdots a_n, b_n\cdots b_1)=
\sum_{k=0}^{n-1} \ff_1(a_1b_{1+k})\cdot\ff_1(a_2b_{2+k})\cdots
\ff_1(a_nb_{n+k}).$$
\end{enumerate}
For a visualization of this formula, one should think of two
concentric circles with the $a$'s on one of them and the $b$'s on
the other. However, whereas on one circle we have a clockwise
orientation of the points, on the other circle the orientation is
counter-clockwise. Thus, in order to match up these points modulo a
rotation of the circles, we have to pair the indices as in the sum
above.

\end{definition}

Recall that in the combinatorial description of freeness \cite{NSp},
the extension of $\ff_1$ to a multiplicative function on
non-crossing partitions plays a fundamental role. In the same way,
second order freeness will rely on a suitable extension of $\ff_2$.

\begin{notation}\label{extension}
Let $(\cA,\ff_1,\ff_2)$ be a second order non-commutative
probability space. Then we extend the definition of $\ff_1$ and
$\ff_2$ as follows:
\begin{align*}
\ff_1:\bigcup_{n=1}^\infty\bigl(S_n\times \cA^n)&\to\CC\\
(\pi,a_1,\dots,a_n)&\mapsto \ff_1(\pi)[a_1,\dots,a_n]
\end{align*}
is, for a cycle $\pi=(i_1, i_2,\dots, i_r)$, given by
$$\ff_1(\pi)[a_1,\dots,a_n]:=\ff_1(a_{i_1}a_{i_2}a_{i_3}\cdots
a_{i_r})$$ and extended to general $\pi\in S_n$ by multiplicativity
$$\ff_1(\pi_1\times\pi_2)[a_1,\dots,a_n]=
\ff_1(\pi_1)[a_1,\dots,a_n]\cdot \ff_1(\pi_2)[a_1,\dots,a_n].$$ In a
similar way,
\begin{align*}
\ff_2:\bigcup_{m,n=1}^\infty\bigl(S_m\times S_n\times
\cA^m\times \cA^n)&\to\CC\\
(\pi_1,\pi_2,a_1,\dots,a_m,b_1,\dots,b_m) &\mapsto
\ff_2(\pi_1,\pi_2)[a_1,\dots,a_m;b_1,\dots,b_n]
\end{align*}
is defined, for two cycles $\pi_1=(i_1, i_2, \dots, i_p)$ and
$\pi_2=(j_1, j_2, \dots,\ab j_r)$, by
$$\ff_2(\pi_1,\pi_2)[a_1,\dots,a_m;b_1,\dots,b_n]:=
\ff_2(a_{i_1}a_{i_2}\cdots a_{i_p},b_{j_1}b_{j_2}\cdots b_{j_r})$$
and extended to the general situation by the derivation
property
\begin{multline}\label{leftderivation}
\ff_2(\pi_1\times\pi_2,\pi_3)[a_1,\dots,a_m;b_1,\dots,b_n]
\\=\ff_2(\pi_1,\pi_3)[a_1,\dots,a_m;b_1,\dots,b_n]\cdot
\ff_1(\pi_2)[a_1,\dots,a_m,b_1,\dots,b_n]
\\+
\ff_2(\pi_2,\pi_3)[a_1,\dots,a_m;b_1,\dots,b_n]\cdot
\ff_1(\pi_1)[a_1,\dots,a_m,b_1,\dots,b_n].
\end{multline}
and
\begin{multline}\label{rightderivation}
\ff_2(\pi_1,\pi_2\times\pi_3)[a_1,\dots,a_m;b_1,\dots,b_n]
\\=\ff_2(\pi_1,\pi_2)[a_1,\dots,a_m;b_1,\dots,b_n]\cdot
\ff_1(\pi_3)[a_1,\dots,a_m,b_1,\dots,b_n]
\\+
\ff_2(\pi_1,\pi_3)[a_1,\dots,a_m;b_1,\dots,b_n]\cdot
\ff_1(\pi_2)[a_1,\dots,a_m,b_1,\dots,b_n].
\end{multline}
\end{notation}

\begin{remark}\label{associativity}
Let $(A_i)_{i \in I}$ be a family of unital subalgebras of
the second order probability space $(A, \ff_1, \ff_2)$ which are
free of second order. Suppose that for each $i$ we have
$(B_{i,j})_{j\in K_i}$ a family of unital subalgebras of $A_i$
which are free of second order. By \cite[Prop. 2.5.5 (iii)]{VDN}
$(B_{i,j})_{j \in \cup_i K_i}$ are free of first order. We leave
as an exercise for the reader to show that the proof of
\cite{VDN} can be adapted to show that  $(B_{i,j})_{j \in \cup_i
K_i}$ are free of second order. \end{remark}

\section{Asymptotic second order freeness for unitary random matrices}

\begin{notation}\label{notation1}
Suppose $\epsilon : [2l] \rightarrow \{-1, 1\}$ is such that
$\sum_{i=1}^{2l} \epsilon_i = 0$. We write $\epsilon^{-1}(1) = \{
p_1, p_2, \dots , p_l\}$ and $\epsilon^{-1}(-1) = \{q_1, q_2, \dots
, q_l\}$, with $p_1 < p_2 < \cdots < p_l$ and $q_1 < q_2 < \cdots <
q_l$. Let $S^{(\epsilon)}_{2l}$ be the permutations $\pi$ in
$S_{2l}$ such that $\pi$ takes $\{p_1,\dots,p_l\}$ onto
$\{q_1,\dots,q_l\}$ and vice versa. Given a $\pi$ in
$S^{(\epsilon)}_{2l}$ we may extract a pair of permutations
$\alpha_\pi$ and $\beta_\pi$ in $S_l$ from the equations 

\begin{equation}\label{pi_definition}
\pi(p_{\alpha_\pi(k)}) = q_{k} \mbox{ and } \pi(q_k) =
p_{\beta_\pi(k)}
\end{equation}
and conversely: $(\alpha, \beta) \mapsto
\pi_{\alpha, \beta}$. Thus we have a bijection of sets between
$S^{(\epsilon)}_{2l}$ and $S_l \times S_l$.

Given $\pi \in S^{(\epsilon)}_{2l}$ we let $\tilde \pi \in S_l$ be
defined by $$ \pi^2(p_k) = p_{\tilde \pi(k)}$$ Note that $\tilde
\pi_{\alpha, \beta} = \beta \alpha^{-1}$.
\end{notation}

Also we have
$$\#(\pi)=\#(\tilde\pi),$$
and thus
$$\vert \pi\vert=\vert\tilde\pi\vert+l.$$

\begin{lemma}
Fix $l\in\NN$ and $\gamma\in S_{2l}$. Let, for $N\in\NN$, $U$ be a
Haar distributed unitary $N\times N$ random matrix. Let $\epsilon :
[2l] \rightarrow \{-1, 1\}$ such that $\sum_{i=1}^{2l} \epsilon_i =
0$. Then we have for all $1\leq
r_1,\dots,r_{2l},s_1,\dots,s_{2l}\leq N$ that
\begin{equation}
\EE\big(U^{\epsilon_1}_{r_1, s_{\gamma(1)}} \cdots\
         U^{\epsilon_{2l}}_{r_{2l}, s_{\gamma(2l)}}
\big) = \sum_{\pi \in S^{(\epsilon)}_{2l}} \prod_{k=1}^{2l}
\delta_{r_k, s_{\gamma(\pi(k))}}\,\Wg(N,\tilde \pi).
\end{equation}
\end{lemma}

\begin{proof}
Let $i_k,i'_k, j_k, j'_k$ ($1 \leq k \leq l$) be such that
$$\EE\big( U^{\epsilon_1}_{r_1, s_{\gamma(1)}} \cdots\
       U^{\epsilon_{2l}}_{r_{2l}, s_{\gamma(2l)}}
\big)
 =
\EE\big(
 U_{i'_1, j'_1} \cdots\ U_{i'_l, j'_l}
       U^{-1}_{j_1, i_1}\cdots\ U^{-1}_{j_l,i_l}
\bigr),
$$
i.e. let $\epsilon^{-1}(1) = \{p_1, \dots , p_l\}$ with $p_1 <
\cdots < p_l$ and $\epsilon^{-1}(-1) = \{q_1 , \dots, q_l \}$
with $q_1 < \cdots < q_l$ and $i'_k = r_{p_k}$, $j'_k =
s_{\gamma(p_k)}$, $i_k = s_{\gamma(q_k)}$, and $j_k =
r_{q_k}$. 

Now suppose that $\alpha$ and $\beta$ in $S_l$ and $\pi \in
S_{2l}^{(\epsilon)}$ is as in equation (\ref{pi_definition})
above.

Thus we have
$$i_k = s_{\gamma(q_k)} =
s_{\gamma(\pi(p_{\alpha(k)}))},\quad\text{and}\quad
i'_{\alpha(k)} = r_{p_{\alpha(k)}},$$ and
$$j'_{\beta(k)}
= s_{\gamma(p_{\beta(k)})} = s_{\gamma(\pi(q_k))},
\quad\text{and}\quad j_k = r_{q_k}$$ which shows that
$$ i_k = i'_{\alpha(k)} \iff r_{p_{\alpha(k)}} =
s_{ \gamma(\pi(p_{\alpha(k)}))  }$$ and
$$j_k = j'_{\beta(k)} \iff r_{q_k} = s_{\gamma(\pi(q_k))}. $$
Thus
$$ \prod_{k=1}^l \delta_{i_k, i'_{\alpha(k)}} \delta_{j_k,
j'_{\beta(k)}} = \prod_{k=1}^{2l} \delta_{r_k,
s_{\gamma(\pi(k))}}.
$$
Hence
\begin{align*}
\EE\big(&
       U^{\epsilon_1}_{r_1, s_{\gamma(1)}} \cdots \
       U^{\epsilon_{2l}}_{r_{2l}, s_{\gamma(2l)}}
\big) =\EE\big(
 U_{i'_1, j'_1} \cdots\ U_{i'_l, j'_l}
       U^{-1}_{j_1, i_1}\cdots\ U^{-1}_{j_l,i_l}
\big)\\
&=\sum_{\alpha,\beta\in S_n} \delta_{i_1i'_{\alpha(1)}}\cdots
\delta_{i_ni'_{\alpha(n)}}
\delta_{j_1j'_{\beta(1)}}\cdots\delta_{j_nj'_{\beta(n)}}
\Wg(\beta\alpha^{-1})\\
&=
 \sum_{\pi \in S^{(\epsilon)}_{2l}}
  \prod_{k=1}^{2l} \delta_{r_k, s_{\gamma(\pi(k))}}\,
                   \Wg(\tilde \pi).
\end{align*}
\end{proof}

We can now address the question how to calculate expectations of
products of traces of our matrices. The following result is exact
for each $N$; later on we will look on its asymptotic version.

Note that the notation $\Tr_\pi(D_1,\dots,D_n)$ for $\pi\in
S_n$ is defined in the usual multiplicative way, as was done
in Notation \ref{extension} for $\ff_1$.

We shall need  the following standard lemma. For $D_1, \dots
,\ab D_p \in M_N(\mathcal{A})$ let the entries of $D_i$ be
$(D^{(i)}_{r,s})$.

\begin{lemma}
Let $\pi \in S_n$ and $D_1, \dots , D_n \in
M_N(\mathcal{A})$. Then 
$$
\Tr_\pi(D_1, D_2, \dots, D_n) = 
\sum_{j_1, j_2, \dots, j_n}
D^{(1)}_{j_1, j_{\pi(1)}} 
D^{(2)}_{j_2, j_{\pi(2)}} \cdots
D^{(n)}_{j_n, j_{\pi(n)}} 
$$
\end{lemma}

Given $m_1, \dots , m_k$, let $\gamma_{m_1, \dots , m_k}$ be
the permutation of $[m_1 + \dots  + m_k]$ with $k$ cycles
where the $i^{th}$ cycle is $(m_1 + \cdots + m_{i-1}+1, \dots
, m_1 + \cdots + m_i)$.

\begin{proposition}
\label{prop:computationofmoments} Fix $m_1,\dots,m_k\in\NN$ such
that $m_1+\cdots+m_k=2l$ is even. Let, for fixed $N\in \NN$, $U$ be
a Haar distributed unitary $N\times N$-random matrix and
$D_1,\dots,D_{2l}$ be $N\times N$-random matrices which are
independent from $U$. Let $\epsilon: [2l] \rightarrow \{-1, 1\}$
with $\sum_{i=1}^{2l} \epsilon_i = 0$. Put
$\gamma=\gamma_{m_1,\dots,m_k}$. Then
\begin{multline}\label{eq:computationofmoments} 
\EE\big(\Tr(D_{1}U^{\epsilon_1}
\cdots D_{m_1}U^{\epsilon_{m_1}}) \\
\times    \Tr(D_{m_1+1}U^{\epsilon_{m_1+1}} \cdots
    D_{m_1+m_2}U^{\epsilon_{m_1+m_2}}) \times \cdots \\
\times  \Tr(D_{m_1+ \cdots + m_{k-1} + 1}U^{\epsilon_{m_1+  \cdots + m_{k-1} + 1}} \cdots
    D_{m_1+ \cdots + m_k}U^{\epsilon_{m_1+ \cdots +
m_k}})\big)\\ 
= \sum_{\pi \in S^{(\epsilon)}_{2l}}
 \Wg(N,\tilde \pi)\cdot \EE\big(\Tr_{\gamma
\pi^{-1}}(D_1,\dots,D_{2l})
\big).
\end{multline}
\end{proposition}

\begin{proof}
Summations over $r$'s and $s$'s in the following formulas are
from 1 to $N$.

\begin{eqnarray*}
\lefteqn{ \EE\big(\Tr(D_{1}U^{\epsilon_1} \cdots
D_{m_1}U^{\epsilon_{m_1}}) 
        \Tr(D_{m_1+1}U^{\epsilon_{m_1+1}} \cdots
    D_{m_1+m_2}U^{\epsilon_{m_1+m_2}}) }\\
&& \cdots \times  
\Tr(D_{m_1+ \cdots + m_{k-1} + 1}U^{\epsilon_{m_1+  
\cdots + m_{k-1} + 1}} \cdots
    D_{m_1+ \cdots + m_k}U^{\epsilon_{m_1+ \cdots +
m_k}})\big)  \\
& = &
\mathop{\sum_{r_1, \dots , r_{2l}}}_{s_1,\dots,s_{2l}} 
         \EE\big(
         U^{\epsilon_1}_{s_1, r_{\gamma(1)}} \cdots\
         U^{\epsilon_{2l}}_{s_{2l}, r_{\gamma(2l)}}
\big)
\cdot \EE\big(D^{(1)}_{r_1, s_1}\cdots D^{(2l)}_{r_{2l}, 
s_{2l}}\big)
 \\
& = &
\mathop{\sum_{r_1, \dots , r_{2l}}}_{s_1,\dots, s_{2l}}\
 \sum_{\pi \in S^{(\epsilon)}_{2l}} \
  \prod_{k=1}^{2l} \delta_{s_k, r_{\gamma(\pi(k))}}\
                   \Wg(\tilde \pi) \cdot
\EE\big(D^{(1)}_{r_1, s_1}\cdots D^{(2l)}_{r_{2l},
s_{2l}}\big)  \\
& = & 
\sum_{\pi \in S^{(\epsilon)}_{2l}}\Wg(\tilde \pi)
 \mathop{\sum_{r_1,\dots, r_{2l}}}_{s_1, \dots , s_{2l}}
  \prod_{k=1}^{2l} \delta_{s_k, r_{\gamma(\pi(k))}}\cdot
\EE\big( D^{(1)}_{r_1, s_1}\cdots D^{(2l)}_{r_{2l}, s_{2l}}
\big) \\
& = & 
\sum_{\pi \in S^{(\epsilon)}_{2l}}\Wg(\tilde \pi)
\sum_{r_1,\dots, r_{2l}}
\EE\big( D^{(1)}_{r_1, r_{\gamma(\pi(1))}}
\cdots D^{(2l)}_{r_{2l}, r_{\gamma(\pi(2l))}}
\big) \\
& = & 
\sum_{\pi \in S^{(\epsilon)}_{2l}} \Wg(\tilde \pi)
\EE\big( \Tr_{\gamma\pi} (D_1,\dots,D_{2l}) \big) \\
& = & 
\sum_{\pi \in S^{(\epsilon)}_{2l}} \Wg(\tilde \pi)
\EE\big( \Tr_{\gamma\pi^{-1}} (D_1,\dots,D_{2l}) \big)
\end{eqnarray*}
In the last equality we used that $\Wg(\tilde\pi)$ depends
only on the conjugacy class of $\pi$.
\end{proof}

Motivated by the result of Voiculescu \cite{Voi1,Voi2} that Haar
distributed unitary random matrices and constant matrices are
asymptotically free, we want to investigate now the corresponding
question for second order freeness. It will turn out that one can
replace the constant matrices by another ensemble of random
matrices, as long as those are independent from the unitary random
matrices. Of course, we have to assume that the second ensemble has
some asymptotic limit distribution. This is formalized in the
following definition. Note that we make a quite strong requirement
on the vanishing of the higher order cumulants. This is however in
accordance with the observation that in many cases the unnormalized
traces converge to Gaussian random variables. Of course, if we have
a non-probabilistic ensemble of constant matrices, then the only
requirement is the convergence of $\kk_1$; all other cumulants are
automatically zero.

\begin{definition} \label{def:asymptoticfreeness}
1) Let $\{A_1,\dots,A_s\}_N$ be a sequence of $N\times N$-random
matrices. We say that they have a \emph{second order limit
distribution} if there exists a second order non-commutative
probability space $(\cA,\ff_1, \ab \ff_2)$ and $a_1,\dots,a_s\in\cA$
such that for all polynomials $p_1,p_2,\dots$ in $s$ non-commuting
indeterminates we have
\begin{equation}\label{asymp1}
\lim_{N\to\infty}\kk_1\big(\tr(p_1(A_1,\dots, A_s)\big)\big)=
\ff_1\bigl(p_1(a_1,\dots, a_s)\bigr),
\end{equation}

\begin{multline}\label{asymp2}
\lim_{N\to\infty} \kk_2\big(\Tr(p_1(A_1,\dots,A_s)),
\Tr(p_2(A_{1},\dots, A_{s}))\big) \\
\mbox{} =
\ff_2\bigl(p_1(a_1,\dots,a_s), p_2 (a_{1},\dots, a_{s})\bigr),
\end{multline}
and, for $r\geq 3$,
\begin{equation}\label{asymp3}
\lim_{N\to\infty} \kk_r\big( \Tr( p_1(A_1,\dots, A_s)),\dots,
\Tr( p_r(A_1,\dots, A_s)) \big) = 0.
\end{equation}
2) We say that two sequences of $N\times N$-random matrices,
$\{A_1,\dots,A_s\}_N$ and $\{B_1,\dots,B_t\}_N$, are
\emph{asymptotically free of second order} if the sequence
$\{A_1,\dots,A_s,B_1,\dots,B_t\}_N$ has a second order limit
distribution, given by $(\cA,\ff_1,\ff_2)$ and $a_1,\dots,a_s,
b_1,\dots,b_t\in\cA$, and if the unital algebras
$$\cA_1:=\alg(1,a_1,\dots,a_s)\qquad\text{and}\qquad
\cA_2:=\alg(1,b_1,\dots,b_t)$$ are free with respect to
$(\ff_1,\ff_2)$.
\end{definition}

\begin{notation}
Fix $m,n\in\NN$ and let $\epsilon:[1,m+n]\to \{-1,+1\}$. We defined
$S_{m+n}^\he$ in Notation \ref{notation1}, for the case where
$\sum_{k=1}^{m+n} \epsilon(k)=0$, as those permutations in $S_{m+n}$
for which $\epsilon$ alternates cyclically between $-1$ and $+1$ on
all cycles. Note that this definition also makes sense in the case
where the sum of the $\epsilon$'s is not equal to zero, then we just
have $S_{m+n}^\he=\emptyset$. Let $\epsilon_1$ and $\epsilon_2$ be
the restrictions of $\epsilon$ to $[1,m]$ and to $[m+1,m+n]$,
respectively. Then we put
$$\SNCe(m,n):=S_{m+n}^\he\cap \SNC(m,n)$$
and
$$NC^{(\epsilon_1)}(m):=S_m^{(\epsilon_1)}\cap \NC(m),\qquad
NC^{(\epsilon_2)}(n):=S_n^{(\epsilon_2)}\cap \NC(n).$$
\end{notation}

\begin{theorem}\label{main}
Let $\{U\}_N$ be a sequence of Haar distributed unitary $N\times
N$-random matrices and $\{A_1,\dots,A_s\}_N$ a sequence of $N\times
N$-random matrices which has a second order limit distribution,
given by $(\cA,\ff_1,\ff_2)$ and $a_1,\dots,a_s\in\cA$. Furthermore,
assume that $\{U\}_N$ and $\{A_1,\dots,A_s\}_N$ are independent. Fix
now $m,n\in\NN$ and consider polynomials $p_1,\dots,p_{m+n}$ in $s$
non-commuting indeterminates. If we put $($for $i=1,\dots,m+n)$
$$D_i:=p_i(A_1,\dots,A_s)\qquad\text{and}\qquad
d_i:=p_i(a_1,\dots,a_s),$$ then we have for all
$\epsilon(1),\dots,\epsilon(m+n)\in\{-1,+1\}$ that
\begin{align}\label{main:eq}
&\lim_{N\to\infty} \kk_2\big(
\Tr(D_{1}U^{\epsilon_1} \cdots
D_{m}U^{\epsilon_m}),
        \Tr(D_{m+1}U^{\epsilon_{m+1}} \cdots
    D_{m+n}U^{\epsilon_{m+n}})\big) \\ \notag
&= \sum_{\pi\in\SNCe(m,n)} \mu(\tilde \pi)\cdot
\ff_1(\gamma_{m,n}\pi^{-1})[d_1,\dots,d_{m+n}]    \\ \notag
&\quad+
\mathop{\sum_{\pi_1\in NC^{(\epsilon_1)}(m)}}_{\pi_2\in
NC^{(\epsilon_2)}(n)} 
\Bigl(\mu_2(\tilde \pi_1, \tilde \pi_2)\cdot
\ff_1(\gamma_m\pi_1^{-1}\times \gamma_n\pi_2^{-1})[d_1,\dots,
d_{m+n}]\\
\notag 
&  \quad\qquad\qquad\qquad+
\mu(\tilde \pi_1 \times\tilde \pi_2)\cdot
\ff_2(\gamma_m\pi_1^{-1},\gamma_n\pi_2^{-1})
[d_1,\dots, d_{m+n}]\Bigr).
\end{align}
\end{theorem}

Note that in the case where the sum of the $\epsilon$'s is different
from zero this just states that the limit of $\kk_2$ vanishes.

\begin{proof}
For notational convenience, we will sometimes write $m+n=2l$ in the
following, and also use $\gamma:=\gamma_{m,n}$.

We have
\begin{eqnarray*}\lefteqn{
\kk_2\big(\Tr(D_1U^{\epsilon_1} \cdots D_{m}U^{\epsilon_m}),
        \Tr(D_{m+1}U^{\epsilon_{m+1}} \cdots
    D_{2l}U^{\epsilon_{2l}})\big)  }  \\  
& = &
\EE\big( \Tr(D_{1}U^{\epsilon_1} \cdots D_{m}U^{\epsilon_m})
        \Tr(D_{m+1}U^{\epsilon_{m+1}} \cdots
    D_{2l}U^{\epsilon_{2l}})\big)   \\
&& \mbox{} - \EE\big( \Tr(D_{1}U^{\epsilon_1} \cdots
D_{m}U^{\epsilon_m})\big)  \cdot \EE\big(
        \Tr(D_{m+1}U^{\epsilon_{m+1}} \cdots
    D_{2l}U^{\epsilon_{2l}})\big)     \\ 
& = &
\sum_{\pi \in S^{(\epsilon)}_{2l}}
\Wg(\tilde \pi)\cdot 
\EE\big( \Tr_{\gamma \pi^{-1}
}(D_1,\dots,D_{2l})\big) \\ 
&& \mbox{} - 
\mathop{\sum_{\pi_1\in S^{(\epsilon_1)}_{m}}}_{\pi_2\in
S^{(\epsilon_2)}_n} 
\Wg(\tilde \pi_1)\Wg(\tilde \pi_2) \cdot
\EE\big( \Tr_{\gamma_m \pi_1^{-1}}(D_1,\dots,D_{m})\big)
 \\[-10pt] 
&& \mbox{}\hskip10em\cdot
\EE\big( \Tr_{\gamma_n \pi_2^{-1}}
(D_{m+1},\dots,D_{2l}) \big) \\
& = &
\mathop{\sum_{\pi \in 
S^{(\epsilon)}_{2l}}}_{\pi\text{\ connected}}
\Wg(\tilde \pi)\cdot \EE\big(\Tr_{\gamma\pi^{-1}} 
(D_1,\dots,D_{2l})  \big)\\ 
&& \mbox{} + 
\mathop{\sum_{\pi_1\in S^{(\epsilon_1)}_{m}}}_{\pi_2\in
S^{(\epsilon_2)}_n} 
\Bigl(\Wg(\tilde \pi_1\times\tilde \pi_2) \cdot
\EE\big( \Tr_{\gamma_m \pi_1^{-1} \times \gamma_n \pi_2^{-1}}
(D_1,\dots,D_{2l})  \big) \\ 
&& \mbox{} -
\Wg(\tilde \pi_1)\Wg(\tilde \pi_2) 
\EE\big( \Tr_{\gamma_m \pi_1^{-1}}
(D_1,\dots,D_{m}) \big) \\ && \hskip10em\mbox{} \cdot
\EE\big( \Tr_{\gamma_n \pi_2^{-1}}(D_{m+1},\dots,D_{2l})
\big) \Bigr) 
\end{eqnarray*} 

Note that if either $m$ or $n$ is odd then the last two terms
are zero, which is consistent with equation (\ref{main:eq}), as
in this case $NC^{(\epsilon_1)}(m)$ and $NC^{(\epsilon_2)}(n)$
are empty. So for the remainder of the proof we shall assume
that $m$ and $n$ are even.

The leading order in the first summand for a connected $\pi$
is given by
\begin{multline*}
\mu(\tilde\pi)N^{-(\vert\tilde\pi\vert+(m+n)/2)}\cdot
N^{\#(\gamma\pi^{-1})}\cdot \EE\big( \tr_{\gamma
\pi^{-1} }(D_1,\dots,D_{m+n})
\big) = \\
=N^{m+n-\vert\pi\vert-\vert\gamma\pi^{-1}\vert}\cdot
\mu(\tilde\pi)\cdot \EE\big(\tr_{\gamma
\pi^{-1}}( D_1,\dots,D_{m+n}) \big).
\end{multline*}

\noindent
Recall that, for a connected $\pi$, we always have
$$m+n-\vert\pi\vert-\vert\gamma\pi^{-1}\vert\leq 0,$$
and equality is exactly achieved in the case where $\pi$
is annular non-crossing. Thus, in the limit $N\to\infty$ the
first sum gives the contribution
$$
\sum_{\pi\in\SNCe(m,n)}\mu(\tilde \pi)
\cdot\ff_1(\gamma\pi^{-1}) [d_1,\dots,d_{m+n}].
$$

We can rewrite the second sum as
\begin{eqnarray}\label{secondsum}\lefteqn{ 
\mathop{\sum_{\pi_1\in S^{(\epsilon_1)}_{m}}}_{\pi_2\in
S^{(\epsilon_2)}_n} 
\Big\{ \Wg(\tilde \pi_1\times\tilde \pi_2) -
\Wg(\tilde \pi_1)\Wg(\tilde \pi_2) \Big\} } \\[-20pt]\notag
&& \hskip6em\mbox{} \times
\EE\big( \Tr_{\gamma_m \pi_1^{-1} \times \gamma_n \pi_2^{-1}}
(D_1,\dots,D_{m+n})  \big)  \\ \notag
&& \mbox{} +
\Wg(\tilde \pi_1)\Wg(\tilde \pi_2) \Big\{  
\EE\big( \Tr_{\gamma_m \pi_1^{-1} \times \gamma_n \pi_2^{-1}}
(D_1,\dots,D_{m+n})  \big) \\ \notag && \mbox{} -
\EE\big( \Tr_{\gamma_m \pi_1^{-1}}
(D_1,\dots,D_{m}) \big) \cdot
\EE\big( \Tr_{\gamma_n \pi_2^{-1}}(D_{m+1},\dots,D_{m+n})
\big)  \Big\}
\end{eqnarray}

For a disconnected $\pi_1\times\pi_2$ the
leading orders in $N$ of all relevant terms are given as
follows: 
$\Wg(\tilde \pi_1\times\tilde \pi_2)$ and
$\Wg(\tilde\pi_1)\Wg(\tilde\pi_2)$ both have leading order
(note that $\mu$ is multiplicative)
$$\mu(\tilde\pi_1)\mu(\tilde\pi_2)N^{-(m+n) +\#(\pi_1)
+ \#(\pi_2)};$$

\noindent
$\EE\big( \Tr_{\gamma_m \pi_1^{-1}}( D_1,\dots,D_{m})\big)
\cdot
\EE\big( \Tr_{\gamma_n \pi_2^{-1}}( D_{m+1},\dots,D_{m+n})
\big)$ and \\
$\EE\big(\Tr_{\gamma_m \pi_1^{-1}\times \gamma_n \pi_2^{-1}}
( D_1,\dots,D_{m+n}) \big) $ are both asymptotic to
$$\ff_1(\gamma_m\pi_1^{-1}\times \gamma_n\pi_2^{-1})[d_1,\dots,
d_{m+n}] N^{\#(\gamma_m\pi_1^{-1}) + \#(
\gamma_n\pi_2^{-1})};$$

$\Wg(\tilde \pi_1\times\tilde
\pi_2)-\Wg(\tilde\pi_1)\Wg(\tilde\pi_2)$ has leading order
$$\mu_2(\tilde\pi_1,\tilde\pi_2)
\cdot N^{-(m+n) + \#(\pi_1) + \#(\pi_2) - 2}.
$$

Now 
\begin{multline*}
-(m+n) + \#(\pi_1) + \#(\pi_2) - 2 + \#(\gamma_m
\pi_1^{-1}) + \#( \gamma_n \pi_2^{-1} )\\
= -( m + 1 - \#(\pi_1) - \#(\gamma_m \pi_1^{-1} )) 
-( n + 1 - \#(\pi_2) - \#(\gamma_n \pi_2^{-1} )) \leq 0
\end{multline*}
with equality only if both $\pi_1 \in NC^{(\epsilon_1)}(m)$
and $\pi_2 \in NC^{(\epsilon_2)}(n)$.

Thus
$$
\lim_N \Big\{ \Wg(\tilde \pi_1 \times \tilde \pi_2) -
\Wg(\tilde  \pi_1) \Wg(\tilde \pi_2) \Big\} \EE(
\Tr_{\gamma_{m,n}(\pi_1 \times \pi_2)^{-1}} (
D_1, \dots , D_{m+n} )  
$$
$$
 = 
\begin{cases}
\mu_2( \tilde\pi_1, \tilde\pi_2 ) 
\ff_1(\gamma_m\pi_1^{-1}\times \gamma_n\pi_2^{-1})[d_1,\dots,
d_{m+n}] & \begin{cases}
\pi_1 \in NC^{(\epsilon_1)}(m) \\
\textrm{\ and\ } &\\
\pi_2 \in NC^{(\epsilon_2)}(n) 
\end{cases} \\
0 & \textrm{otherwise} \\
\end{cases}
$$

To deal with the second term of the second sum
(\ref{secondsum}) we will use the following notation. Let the
cycles of $\gamma_m \pi_1^{-1}$ be $c_1  \cdots 
c_r$ and the cycles of $\gamma_n\pi_2^{-1}$ be $c_{r+1} 
\cdots  c_{r+s}$. Let $a_i = \Tr_{c_i}(D_1, \dots D_m)$ for $1
\leq i \leq r$ and $b_j = \Tr_{c_{r+j}}(D_{m+1}, \dots ,
D_{m+n})$ for $1 \leq j \leq s$. Then
\begin{multline*}
\kk_2( \Tr_{\gamma_m\pi_1^{-1}}(D_1,     \dots, D_m), 
       \Tr_{\gamma_n\pi_2^{-1}}(D_{m+1}, \dots, D_{m+n}) \\
=
\kk_2(a_1 \cdots a_r, b_1 \cdots b_s)
\end{multline*}

So let us find for which $\pi_1 \in S_m^{(\epsilon_1)}$,
$\pi_2 \in S_n^{(\epsilon_2)}$, and  $\tau \in \cP(r+s)$ we have a
non-zero limit of
\begin{equation}\label{generalterm}
\Wg(\tilde\pi_1)\, \Wg(\tilde\pi_2) \,
\kk_\tau ( a_1, \dots, a_r, b_1, \dots , b_s)
\end{equation}

As noted above the order of $\Wg(\tilde\pi_1)\,
\Wg(\tilde\pi_2)$ is $N^{-(m+n) + \#(\pi_1) + \#(\pi_2)}$. 
By equation (\ref{leonov}) and the definition of a second order limit
distribution, $\kk_\tau ( a_1, \dots, a_r, b_1, \dots , b_s)$ is
$O(N^c)$ where $c$ is the number of singletons of $\tau$. Thus the order of
(\ref{generalterm}) is 
\begin{eqnarray*}\lefteqn{
-(m+n) + \#(\pi_1) + \#(\pi_2) + c }\\
&=& - \big( m + 1 - \#(\pi_1) - \#(\gamma_m\pi_1^{-1})\big) -
 \big(n + 1 - \#(\pi_2) - \#(\gamma_n\pi_2^{-1}) \big) \\
&&\quad \mbox{} + c + 2 - (r + s)
\end{eqnarray*}
Hence (\ref{generalterm}) will vanish unless three conditions
are satisfied: we must have that $\pi_1$ and $\pi_2$ are
non-crossing and $c = r + s -2$, i.e. $\tau$ has one pair and
the rest of its blocks are singletons. 

Thus
\begin{eqnarray*}\label{secondsumlimit}
\lefteqn{ \lim_N
\sum_{\pi_1\in S^{(\epsilon_1)}_{m}, \pi_2\in
S^{(\epsilon_2)}_n} 
\Wg(\tilde \pi_1)\Wg(\tilde \pi_2) } \\ 
&& \qquad \mbox{} \times
\kk_2 \Big( \Tr_{\gamma_m \pi_1^{-1}}
            (D_1,\dots,D_{m}) ,
      \Tr_{\gamma_n \pi_2^{-1}}
            (D_{m+1},\dots,D_{2l})
  \Big) \\
&=&
\mathop{\sum_{\pi_1\in NC^{(\epsilon_1)}(m)}}_{\pi_2\in
NC^{(\epsilon_2)}(n)} \mu(\tilde\pi_1) \mu(\tilde\pi_2) \\
&& \qquad  \mbox{} \times
\lim_N N^{r+s-2} \sum_{\tau \in \cP(r+s)}
\kk_\tau (a_1, \dots , a_r, b_1, \dots , b_s)
\end{eqnarray*} 
where the $\tau$'s in the sum have one pair and the remainder are
singletons and $\tau \vee 1_{r,s} = 1_{r+s}$.

So the remainder of the proof is to show that
\begin{multline*}
\lim_N  N^{r+s-2} \mathop{\sum_{\tau \in \cP(r+s)}}
\kk_\tau (a_1, \dots , a_r, b_1, \dots , b_s)\\
=
\ff_2(\gamma_m\pi_1^{-1},\gamma_n\pi_2^{-1})
[d_1,\dots, d_{m+n}]
\end{multline*}
where the sum runs over $\tau$'s as above.

Let $\tau \in \cP(r+s)$ be as above with pair $(i,j)$ where
$1 \leq i \leq r$ and $1 \leq j \leq s$ and all other blocks
singletons. Then 
\begin{eqnarray*}\lefteqn{
\kk_\tau( a_1, \dots , a_r, b_1, \dots , b_s) }\\
&=&
\kk_1(a_1) \cdots \widehat{\kk_1(a_i)} \cdots \kk_1(a_r)
\kk_1(b_1) \cdots \widehat{\kk_1(b_j)} \cdots \kk_1(b_s)
\kk_2(a_i, b_j)
\end{eqnarray*}
where the hatted elements are deleted. So after multiplying by
$N^{r+s-2}$
and taking a limit we get
(after omitting the arguments $d_1, \dots , d_{m+n}$ which are
the same for each factor)
$$
\ff_1(c_1) \cdots \widehat{\ff_1(c_i)} \cdots \ff_1(c_r)
\ff_1(c_{r+1}) \cdots \widehat{\ff_1(c_{r+j})} \cdots
\ff_1 (c_{r + s}) \ff_2(c_i, c_{r+j})
$$
Now summing over all $\tau$, which is equivalent to summing
over all $i$ and $j$, we get via the derivation property of
$\ff_2$ (see equations (\ref{leftderivation}) and
(\ref{rightderivation}))
$$
\ff_2(\gamma_m \pi_1^{-1}, \gamma_n \pi_2^{-1})
[d_1, \dots , d_{m+n}]
$$
as required. \end{proof}

\begin{remark}\label{remark:ds}
When all the $D$'s are equal to 1, equation (\ref{main:eq}) implies
the following well known result of Diaconis and Shahshahani
\cite{DS}: for integers $r$ and $s$ 
\begin{equation}\label{diaconis-shashahani}
\lim_N \kk_2\big( \Tr(U^r), \Tr(U^s) \big) =
\begin{cases}
0   & r \not= - s \\
|r| & r = - s
\end{cases}
\end{equation}
Indeed let $m = |r|$ and for $1 \leq i \leq m$, let $\epsilon_i =
\textrm{sgn}(r)$, where $\textrm{sgn}(r)$ denotes the sign of $r$;
let $n = |s|$ and for $m+1 \leq i \leq m+n$, let $\epsilon_i =
\textrm{sgn}(s)$. Then $\epsilon_1 + \cdots + \epsilon_{m+n} = r
+ s$. So if $r+s \not = 0$ then equation (\ref{main:eq}) says
that  \[
\lim_n \kk_2 \big( \Tr(U^r), \Tr(U^s) \big) = 0
\]
Suppose that $r + s = 0$. The second term on the right hand side of
(\ref{main:eq}) is zero since both $NC^{(\epsilon_1)}(m)$ and
$NC^{(\epsilon_2)}(n)$ are empty. For the first term in
(\ref{main:eq}), note that the only elements of
$S_{NC}^{(\epsilon)}(m,n)$ which connect in this alternating way
are pairings, where each block must contain one $U$ and one
$U^\ast$. This forces $m$ and $n$ to be equal. In that case, we
have the freedom of pairing the first $U$ with any of the $n$
$U^\ast$'s. After this choice is made, the rest is determined. Thus
there are $n$ possibilities for such pairings. Since $\mu(\tilde
\pi)$ is always 1 for a pairing we get the claimed formula. \qed
\end{remark}

Let $\epsilon: [2l] \rightarrow \{-1, 1\}$ be such that
$\sum_{i = 1}^{2l} \epsilon_i = 0$. For $\pi \in
S^{(\epsilon)}_{2l}$ let $\tilde\pi \in S_l$ be as in 3.1. A
$\pi$-invariant partition $A$ of $[2l]$ gives $\tilde A$, a
$\tilde\pi$-invariant  partition of $[l]$ as follows. For each
block $V$ of $A$ let $\tilde V = \{ k \mid p_k \in V\}$, where
we have used the notation of 3.1. Also each
$\tilde\pi$-invariant partition of $[l]$ comes from a unique
$\pi$-invariant partition of $[2l]$.

Let $\moeb$ be the M\"obius function on the partially ordered
set of partitions of $[l]$ ordered by inclusion. Let $\pi \in
S_l$ and $A$ be a $\pi$-invariant partition of $[l]$. In
\cite[\S 2.3]{Col} Collins denotes the relative cumulant by
$C_{\Pi_\pi, A}(\pi,N)$, which we will denote by $C_{\pi,A}$. In our
notation  $$
C_{\pi, A} =
\mathop{\sum_{C \in [ \pi, A]}}_{C = \{V_1, \dots , V_k\} }
\moeb(C, A) \Wg(\pi|_{V_1}) \cdots
\Wg(\pi|_{V_k})
$$
where $\pi|_{V_i}$ denotes the restriction of $\pi$ to the
invariant subset $V_i$ and where necessary we have identified $\pi$ with the
partition given by its cycles. Conversely given $A = \{V_1, \dots , V_k \}$ a
$\pi$-invariant partition of $[l]$ we write $\Wg_A(\pi)$ for $\Wg(\pi|_{V_1})
\cdots \Wg(\pi|_{V_k})$. Then by M\"obius inversion we have
$$
\Wg_A(\pi) =
\sum_{C \in [\pi, A]} C_{\pi, C}
$$

\begin{remark} 
When $\pi \in S_{2l}^{(\epsilon)}$ and $A \in \cP(2l)$ is
$\pi$-invariant the equation above can also be written
\begin{equation}\label{weingartencumulant}
\Wg_A(\tilde\pi) =
\sum_{C \in [\pi, A]} C_{\tilde\pi, \tilde
C}
\end{equation}
In \cite[Cor. 2.9]{Col} Collins showed that the order of
$C_{\tilde\pi, \tilde C}$ is at most
$N^{-2l - \#(\pi) + 2 \#(C)}$.  \end{remark}

In the following we address the estimates for higher order
cumulants, $\kk_r$ for $r\geq 3$.

If $D_1,\dots,D_{2l}$ are random matrices and $\pi\in S_{2l}$
is a permutation with cycle structure $\pi=\pi_1 \times
\cdots \times \pi_r$ with
$\pi_i=(\pi_{i,1},\dots,\pi_{i,l(i)})$ we denote $$
\kk_{\pi}(D_1,\dots,D_l)= \kk_r \big( \Tr (D_{\pi_{1,1}}
\cdots D_{\pi_{1,l(1)}}), \Tr (D_{\pi_{2,1}} \cdots
D_{\pi_{2,l(2)}} ),\dots \big). $$

When $A=\{A_1,\dots,A_k\}$ is a
$\pi$-invariant partition of $[2l]$ we can write
$\pi=\pi_1\times \cdots \times\pi_k$ where $\pi_i
= \pi|_{A_i}$  is a permutation of the set $A_i$. We denote
the multiplicative extension of $\kk_\pi$ by
$$ \kk_{\pi,A}(D_1,\dots,D_{2l})=
\kk_{\pi_1} (D_1,\dots,D_{2l})\cdots \kk_{\pi_k}
(D_{1},\dots,D_{2l}) $$

M\"obius inversion gives us that
$$ \EE \big( \Tr_{\pi} (D_1,\dots,D_{2l}) \big)
= 
\mathop{\sum_{A \in \cP(2l)}}_{A\ \pi\textrm{-inv.}}
\kk_{\pi,A} (D_1,\dots,D_{2l})
$$ 
where
the sums run over all $\pi$-invariant partitions $A$ in
$\cP(2l)$.

\begin{theorem}\label{main:higherorder}
Let $\{U\}_N$ be a sequence of Haar distributed unitary
$N\times N$-random matrices and $\{A_1,\dots,A_s\}_N$ a
sequence of $N\times N$-random matrices which has a second
order limit distribution, given by $(\cA,\ff_1,\ff_2)$ and
$a_1,\dots,a_s\in\cA$. Furthermore, assume that $\{U\}_N$ and
$\{A_1,\dots,A_s\}_N$ are independent.  

Suppose $r > 1$ and $m_1,\cdots,m_r$ are positive integers such that
$m_1+\cdots+m_r = 2l$ and $\epsilon_1, \dots, \epsilon_{2l} \in
\{-1,+1\}$ are such that $\sum_{i=1}^{2l} \epsilon_i = 0$.
Consider polynomials
$p_1,\dots,p_{2l}$ in $s$ non-commuting indeterminates. For
$i=1,\dots, 2l$ we set

$$D_i:=p_i(A_1,\dots,A_s)$$ and for $1 \leq i
\leq r$ let
$$
X_i =
\Tr\big( D_{m_1 + \cdots + m_{i-1} +1} U^{\epsilon(m_1 + \cdots
+ m_{i-1} +1)} \cdots
D_{m_1 + \cdots + m_i} U^{\epsilon(m_1 + \cdots +
m_i)} \big)
$$

Then
\begin{equation}\label{main:eqhigherorder}
\kk_r\big(X_1, \dots , X_r)
= \sum_{\pi\in S^{(\epsilon)}_{n}} 
  \sum_{\substack{A,B \\ A\vee B =1_{[1,2l]}}} \kern-1em
C_{\tilde\pi, \tilde A}
  \cdot
\kk_{\gamma \pi^{-1},B} (D_1,\dots,D_{2l}),
\end{equation}
where the second sum runs over pairs $(A,B)$ of partitions of
$[1,2l]$ such that $A$ is $\pi$-invariant and $B$ is
$\gamma\pi^{-1}$-invariant and furthermore $A\vee
B=1_{[1,2l]}$.

Secondly, we have for $r\geq 3$ that
\begin{equation} \label{main:eqhigherorderB}
\lim_{N\to\infty} \kk_r\big( X_1, \dots, X_r \big) = 0.
\end{equation}
If we have $m_1, \dots , m_r$ for which $m_1 + \cdots + m_r$ is
odd or $\epsilon_1, \dots , \epsilon_{2l}$ for which
$\sum_{i=1}^{2l} \epsilon_i \not = 0$ then $\kk_r(X_1, \dots ,
X_r) = 0$. 
\end{theorem}

\begin{proof}
 In order to simplify the writing we shall write $\vec D$ for
$(D_1, \dots ,\ab D_{2l})$. Let $I_i $ be the interval $[m_1 +
\cdots + m_{i-1} + 1, m_1 + \cdots + m_i]$ and for any subset $V
\subset [r]$, $I_V = \cup_{j \in V} I_j$.

If $C = \{V_1, \dots , V_k\}$ is a
partition of $[r]$ we let $S_{V_i}^{(\epsilon)}$ the set of
permutations of $I_{V_i}$ that take $\epsilon^{-1}(1) \cap
I_{V_i}$ onto $\epsilon^{-1}(-1) \cap I_{V_i}$. If these two
sets have different cardinalities then $S_{V_i}^{(\epsilon)}$
is empty. Let $1_C$ be the
partition $\{ I_{V_1}, \dots , I_{V_k} \}$ of $[2l]$. Let
$\gamma_i$ be the cyclic permutation of $I_i$ given by $(m_1 +
\cdots + m_{i-1} + 1, \dots, m_1 + \cdots + m_i)$

With this notation
\begin{eqnarray*}\lefteqn{
\EE_C( X_1, \dots , X_r) = \EE_{V_1}( X_1, \dots , X_r)
\cdots \EE_{V_k}( X_1, \dots , X_r) }\\
&=& \kern-0.5em
    \sum_{\pi_1 \in S_{V_1}^{(\epsilon)}} 
     \kern-0.5em\cdots\kern-0.5em
    \sum_{\pi_k \in S_{V_k}^{(\epsilon)}}
\Wg(\tilde\pi_1) \cdots \Wg(\tilde\pi_k)
\EE\big( \Tr_{\gamma_1\pi_1^{-1}}( \vec D) \big) \cdots
\EE\big( \Tr_{\gamma_k\pi_k^{-1}}( \vec D) \big) \\
&=& \mathop{\sum_{\pi \in S_{2l}^{(\epsilon)}}}_{1_c
\pi\textrm{-inv.}}
  \Wg_{1_C}(\tilde\pi)
  \EE_C \big( \Tr_{\gamma\pi^{-1}}( \vec D ) \big) \\
&=& \mathop{\sum_{\pi \in S_{2l}^{(\epsilon)}}}_{1_c
\pi\textrm{-inv.}}
\sum_{A \in [\pi, 1_C]} C_{\tilde\pi, \tilde
A}
\mathop{\sum_{B \in \cP(2l)}}_{B\ \gamma\pi^{-1}\textrm{-inv.}}
\kk_{\gamma\pi^{-1},B} (\vec D) 
\end{eqnarray*}
Thus
\begin{eqnarray*}\lefteqn{
\kk_r(X_1, \dots , X_r) }\\
&=& \sum_{C \in \cP(r)} \moeb(C, 1_r) \EE_C(X_1, \dots X_r) \\
&=& \sum_{C \in \cP(r)} \moeb(C, 1_r)
\mathop{\sum_{\pi \in S_{2l}^{(\epsilon)}}}_{1_C
\pi\textrm{-inv.}}
\sum_{A \in [\pi, 1_C]} C_{\tilde\pi, \tilde
A}
\mathop{\sum_{B \in \cP(2l)}}_{B\ \gamma\pi^{-1}\textrm{-inv.}}
\kk_{\gamma\pi^{-1},B} (\vec D) \\
&=& \sum_{\pi \in S_{2l}^{(\epsilon)}}
\mathop{\sum_{C \in \cP(r)}}_{1_C \pi\textrm{-inv.}}  
\sum_{A \in [\pi, 1_C]}
\mathop{\sum_{B \in \cP(2l)}}_{B\ \gamma\pi^{-1}
\textrm{-inv.}} 
\moeb(C, 1_r)\, C_{\tilde\pi, \tilde A}\,
\kk_{\gamma\pi^{-1},B} (\vec D) \\
&=& \sum_{\pi \in S_{2l}^{(\epsilon)}}
\mathop{\sum_{A \in \cP(2l)}}_{\pi\textrm{-inv.}}
\mathop{\sum_{B \in \cP(2l)}}_{B\ \gamma\pi^{-1}\textrm{-inv.}}
\mathop{\sum_{C \in \cP(r)}}_{A, B \leq 1_C} 
\moeb(C, 1_r)\, C_{\tilde\pi, \tilde A}\,
\kk_{\gamma\pi^{-1},B} (\vec D) \\
&=& \sum_{\pi \in S_{2l}^{(\epsilon)}}
\mathop{\sum_{A, B \in \cP(2l)}}_{A \vee B = 1_{2l}}
C_{\tilde\pi, \tilde A}\,
\kk_{\gamma\pi^{-1},B} (\vec D)
\end{eqnarray*}
where the sum is over all $A$ and $B$ which are $\pi$ and
$\gamma\pi^{-1}$-invariant respectively. The last equality
followed from the identity
$$
\mathop{\sum_{C \in \cP(r)}}_{A, B \leq 1_C} 
\moeb(C, 1_r) =
\begin{cases}
1 & A \vee B = 1_{2l} \\
0 & \textrm{otherwise}
\end{cases}
$$
This proves (\ref{main:eqhigherorder}).

We know that the order of $C_{\tilde\pi, \tilde
A}$ is $N^{-2l - \#(\pi) + 2 \#(A)}$. Let $c_i$
be the number of blocks of $B$ that contain $i$ cycles of
$\gamma\pi^{-1}$. By our assumption on the second order
limiting distribution of $\{A_1,\dots,A_s\}_N$
$$
\kk_{\gamma\pi^{-1}, B}(\vec D) =
\begin{cases}
O(N^{c_1})  & c_3 + c_4 + \cdots =0\\
o(N^{c_1})  & c_3 + c_4 + \cdots >0
\end{cases}
$$

Suppose first that $c_3 + c_4 + \cdots >0$. Then 
$$\sum_{i\geq 2} i c_i=
(c_2+c_3+\cdots)+\sum_{i\geq 1} (i-1) c_i\geq 1+
\#(\gamma\pi^{-1})-\#(B)$$
So
$$
c_1 = \#(\gamma\pi^{-1}) - \sum_{i\geq 2} i c_i \leq \#(B) - 1
$$
By Lemma \ref{triangleinequality} (1), $\#(A) + \#(B) \leq
2 l + 1$. Thus
$$
-2l - \#(\pi) + 2 \#(A) + c_1 \leq
-2l - \#(\pi) + 2 \#(A) + \#(B) - 1 \leq 0
$$
Hence $C_{\tilde\pi, \tilde A} \cdot
\kk_{\gamma\pi^{-1}}(\vec D) = o(N^0)$ as required.

So now suppose that $c_3 + c_4 + \cdots =0$. Then
$\#(\gamma\pi^{-1}) = c_2 + \#(B)$. In this case 
$$
C_{\tilde\pi, \tilde A} \cdot
\kk_{\gamma\pi^{-1}}(\vec D)
= O(N^{-2l - \#(\pi) + 2 \#(A) + c_1})
$$
and thus it remains to show that
\begin{equation}\label{thefinalbit}
-2l - \#(\pi) + 2\#(A) + c_1 \leq 2 - r
\end{equation}
Note that
$$
|(1_{2l}, \gamma)| = 2 |1_{2l}| - |\gamma| = 2l -2 + r, 
$$
$$
|(A, \pi)| = 2 |A| - |\pi| = 2l - 2\#(A) + \#(\pi)
$$
and
$$
|(B, \gamma\pi^{-1})| = 2 |B| - |\gamma\pi^{-1}|
= 2l - 2 \#(B) + \#(\gamma\pi^{-1})
$$
So by Lemma \ref{triangleinequality} (2) 
\begin{equation}\label{thetriangleinequalityunfolded}
2\big( \#(A) + \#(B) ) - 2l - \#(\pi) - \#(\gamma\pi^{-1})
\leq 2 - r
\end{equation}
However 
$$
2 \#(B) - \#(\gamma\pi^{-1}) = \#(B) - c_2 = c_1
$$
together with (\ref{thetriangleinequalityunfolded}) proves (\ref{thefinalbit})
\end{proof}

\begin{remark}
As a corollary of Theorem \ref{main:higherorder} we obtain that if
$\{ U \}_N$ is a sequence of Haar distributed unitary random
matrices, then $\{ U \}_N$ has a second order limit distribution
given by equation (\ref{diaconis-shashahani}). Indeed, relative to
$\EE\big( \Tr( \cdot ) \big)$, $U$ is already a Haar unitary so
condition (\ref{asymp1}) of Definition \ref{def:asymptoticfreeness}
is satisfied. We have observed in Remark \ref{remark:ds} that
condition (\ref{asymp2}) is satisfied and by Theorem
\ref{main:higherorder} above we have that condition (\ref{asymp3})
is satisfied. 

Then $\{U\}_N$ has a second order limit
distribution which is given by
\begin{align}
\lim_{N\to\infty} & \kk_2\big(\Tr(U^{\epsilon_1} \cdots
U^{\epsilon_m}),
        \Tr(U^{\epsilon_{m+1}} \cdots
    U^{\epsilon_{m+n}})\big)\\ \notag
&= \sum_{\pi\in\SNCe(m,n)} \mu(\tilde \pi) + 
\mathop{\sum_{\pi_1\in NC^{(\epsilon_1)}(m)}}_{%
\pi_2\in NC^{(\epsilon_2)}(n)}
\mu_2(\tilde\pi_1, \tilde \pi_2)
\end{align}

Combining this formula with equation (\ref{diaconis-shashahani})
allows one to derive the values of $\mu_2$. These kind of
questions will be considered elsewhere.\qed
\end{remark}

\begin{theorem}\label{main:free}
Let $\{U\}_N$ be a sequence of Haar distributed unitary $N\times N$
random matrices and $\{A_1,\dots,A_s\}_N$ a sequence of $N\times
N$ random matrices which has a second order limit distribution. If
$\{U\}_N$ and $\{A_1,\dots,A_s\}_N$ are independent, then they are
asymptotically free of second order.
\end{theorem}

\begin{proof}
The asymptotic freeness with respect to $\kk_1\big( \tr(\cdot)\big)$
is essentially the same argument as Voiculescu's proof
\cite{Voi1,Voi2} for the case of constant matrices, see also the
proof of Collins \cite{Col}.

Theorem \ref{main:higherorder} provides the bound on higher order
cumulants so we need to prove now only the second order statement.

We have to consider cyclically alternating and centered words in the
$U$'s and the $A$'s. For the $U$'s, every centered word is a linear
combination of non-trivial powers of $U$, thus it suffices to
consider such powers. Thus we have to look at expressions of the
form
\begin{equation}\label{express}
\kk_2\big(\Tr(B_1 U^{i(1)}\cdots B_p
U^{i(p)}),\Tr(U^{j(r)}C_r\cdots U^{j(1)}C_1)\big),
\end{equation}
where the $B$'s and the $C$'s are centered polynomials in the $A$'s
and $i(1),\dots,i(p),j(1),\dots,j(r)$ are integers different from
zero. We have to show that in the limit $N\to\infty$ the expression
(\ref{express}) converges to
\begin{equation} \lim_N
\delta_{pr}\sum_{k=0}^{p-1}\ff_1(B_1C_{1+k})\ff_1(U^{i(1)}U^{j(1+k)})
\cdots \ff_1(B_pC_{p+k})\ff_1(U^{i(p)}U^{j(p+k)}).
\end{equation}
We can bring the expression (\ref{express}) into the form considered
in Theorem \ref{main} by inserting $1$'s between neighbouring
factors $U$ or neighbouring factors $U^*$. If we relabel the $B$'s,
$C$'s, and $1$'s as $D$'s then we have to look at the following
situation: For polynomials $p_i$ in $s$ non-commuting indeterminates
we consider
$$D_i:=p_i(A_1,\dots,A_s),$$
which are either asymptotically centered or equal to 1. The latter
case can only appear if we have cyclically the pattern $\dots
UD_iU\dots$ or $\dots U^*D_iU^*\dots$. Formally, this means:
\begin{itemize}
\item if $\epsilon_{\gamma^{-1}(i)}=\epsilon_i$ then either
$D_i=1$ (for all $N$, i.e., $p_i=1$) or
$$\lim_{N\to\infty}\kk_1\big(\tr[D_i]\big)=0.$$
\item if $\epsilon_{\gamma^{-1}(i)}\not=\epsilon_i$ then
$$\lim_{N\to\infty}\kk_1\big(\tr[D_i]\big)=0.$$
\end{itemize}
We can now use Theorem \ref{main} for calculating the limit
$$
\lim_{N\to\infty} \kk_2\big(\Tr(D_{1}U^{\epsilon_1} \cdots
D_{m}U^{\epsilon_m}),
        \Tr(D_{m+1}U^{\epsilon_{m+1}} \cdots
    D_{m+n}U^{\epsilon_{m+n}})\big),
$$
and we will argue that most terms appearing there will vanish.
First consider the last two sums in equation (\ref{main:eq}),
corresponding to $\pi_1\in NC(m)$ and $\pi_2\in NC(n)$. Since
$\pi_1$ is non-crossing we have that
$\#(\pi_1)+\#(\gamma_m\pi_1^{-1})=m+1$. Since each cycle of $\pi_1$
must contain at least one $U$ and one $U^*$, we have
$$\#(\pi_1)\leq \frac m2,$$ which implies
$\#(\gamma_m\pi_1^{-1})\geq m/2+1$. However, this can only be true
if $\gamma_m\pi_1^{-1}$ contains at least two singletons. Note that
if $(i)$ is a singleton of $\gamma_m\pi_1^{-1}$ and if we have
$D_i=1$ for that $i$, then we have
$$\gamma_m\pi_1^{-1}(i)=i,\qquad\text{thus}\qquad
\pi_1^{-1}(i)=\gamma_m^{-1}(i)=\gamma^{-1}(i),$$ and hence
$$\epsilon_{\pi_1^{-1}(i)}=\epsilon_{\gamma^{-1}(i)}=\epsilon_i,$$
which is not allowed because $\pi_1$ is from $NC^{(\epsilon_1)}(m)$,
i.e., it must connect alternatingly $U$ with $U^*$. Hence $D_i \not
= 1$ and so $\ff_1(d_i) = \lim_N \kk_1(\tr(D_i))\ab = 0$.
Thus, both
$$\ff_1(\gamma_m\pi_1^{-1}\times \gamma_n\pi_2^{-1})[d_1,\dots,
d_{m+n}]$$ and
$$\ff_2(\gamma_m\pi_1^{-1}, \gamma_n\pi_2^{-1})[d_1,\dots,
d_{m+n}]$$ are zero, because at least one singleton $(i)$
gives the contribution $\ff_1(d_i) \ab = 0$.

Consider now the first summand of equation (\ref{main:eq}). Suppose
$\pi\in\SNCe(m,n)$. Let us again put $\gamma:=\gamma_{m,n}$.
Since $\pi$ is annular non-crossing we have
$$\vert\pi\vert+\vert\gamma\pi^{-1}\vert=m+n,$$
or
$$\#(\pi) +\#(\gamma\pi^{-1})=m+n.$$
Again, each cycle of $\pi$ must contain at least two elements, i.e.,
$$\#(\pi) \leq \frac{m+n}2,$$
thus
$$\#(\gamma\pi^{-1})\geq\frac{m+n}2.$$
If $\gamma\pi^{-1}$ has a singleton $(i)$, then this will contribute
$\ff_1(d_i)$ and since, as above the case $d_i=1$ is excluded for a
singleton, we get a vanishing contribution in this case. This
implies that, in order to get a non-vanishing contribution,
$\gamma\pi^{-1}$ must contain no singletons, which, however, means
that we must have
$$\#(\gamma\pi^{-1})=\frac{m+n}2,\qquad\text{and thus also}\qquad
\#(\pi)=\frac{m+n}2$$ i.e., all cycles of $\gamma\pi^{-1}$ and of
$\pi$ contain exactly two elements. This, however, can only be the
case if each cycle connects one point on the outer circle to one
point on the inner circle. Being non-crossing fixes the permutation
up to a rotation of the inner circle. Thus, in order to get a
non-vanishing contribution, we need $m=n$ and
$$\pi=(1,\gamma^k(2n))(2,\gamma^k(2n-1)),\dots,(n,\gamma^k(n+1))$$
for some $k=0,1,\dots,n-1$. Note that $\pi$ must always couple a $U$
with a $U^*$ and the factor $\mu(\tilde \pi)$ is always 1 for such
pairings. This gives exactly the contribution as needed for second
order freeness.
\end{proof}

Of course, a natural question in this context is how the 
result of Diaconis and Shahshahani (Remark \ref{remark:ds}) generalizes to the
case of several independent unitary random matrices. Note that as we have
established the existence of a second order limit distribution for Haar
distributed unitary random matrices we can use an independent copy of them as
the ensemble $\{A_1,\dots,A_s\}$ in our Theorem \ref{main:free}. By Remark
\ref{associativity} this can be iterated to give the following.

\begin{theorem}
Let $\{U^{(1)}\}_N,\dots,\{U^{(r)}\}_N$ be $r$ sequences of Haar
distributed unitary $N\times N$-random matrices. If
$\{U^{(1)}\}_N,\dots,\{U^{(r)}\}_N$ are independent, then they are
asymptotically free of second order.
\end{theorem}

This contains the information about the fluctuation of several
independent Haar distributed unitary random matrices. Again, it
suffices to consider traces of reduced words in our random matrices,
i.e., expressions of the form
\begin{equation}\label{words}
\Tr(U_{i(1)}^{k(1)}\cdots U_{i(n)}^{k(n)})
\end{equation}
for $n\in\NN$, and $k(r)\in\ZZ\backslash \{0\}$ and $i(r)\not=
i(r+1)$ for all $r=1,\dots,n$ (where $i(n+1)=i(1)$). But these are
now products in cyclically alternating and centered variables, so
that by the very definition of second order freeness we get
\begin{align}
\lim_{N\to\infty} 
\kk_2\big(&\Tr(U_{i(1)}^{k(1)}\cdots
U_{i(m)}^{k(m)}), \Tr(U_{j(n)}^{l(n)}\cdots
U_{j(1)}^{l(1)})\big) \\ \notag 
&=\delta_{mn} \sum_{r=0}^{n-1}
\ff_1\bigl(U_{i(1)}^{k(1)}U_{j(1+r)}^{l(1+r)}\bigr) \cdots
 \ff_1\bigl(U_{i(n)}^{k(n)}U_{j(n+r)}^{l(n+r)}\bigr).
\end{align}
The contribution of $\ff_1$ in these terms vanishes unless the
matrices and their powers match. Note also that the vanishing of
higher cumulants can be rephrased in a more probabilistic language
by saying that the random variables (\ref{words}) converge to a
Gaussian family.

\begin{corollary}
Let $\{U_{(1)}\}_N,\dots,\{U_{(r)}\}_N$ be independent sequences of
Haar distributed unitary $N\times N$-random matrices. Then, the
collection (\ref{words}) of unnormalized traces in cyclically
reduced words in these random matrices converges to a Gaussian
family of centered random variables whose covariance is given by the
number of matchings between the two reduced words,
\begin{align}
\lim_{N\to\infty} &\kk_2\big(\Tr(U_{i(1)}^{k(1)}\cdots
U_{i(m)}^{k(m)}),
\Tr(U_{j(n)}^{l(n)}\cdots U_{j(1)}^{l(1)})\big) \\
&=\notag \delta_{mn}\cdot\#\bigl\{ r\in\{1,\dots,n\}\mid
i(s)=j(s+r),\\\notag &\qquad\qquad\qquad\qquad\qquad\qquad
k(s)=-l(s+r) \quad\forall s=1,\dots,n \bigl\}
\end{align}
\end{corollary}

This result was also obtained independently in the recent work of
R\u{a}dulescu \cite{Rad} around Connes's embedding problem.

The following theorem gives an easy way to construct families of
random matrices which are asymptotically free of second order.

\begin{theorem}
Let $\{U\}_N$ be a sequence of Haar distributed unitary $N\times
N$-random matrices. Suppose that $\{A_1,\dots,A_s\}_N$ and
$\{B_1,\dots, \ab B_t\}_N$ are sequences of $N \times N$ random matrices each of
which has a second order limit distribution.  Furthermore, assume that $\{A_1,
\ab \dots, \ab A_s\}_N$, $\{B_1,\dots,B_t\}_N$, and $\{U\}_N$ are independent.
Then the sequences $\{A_1,\dots, \ab A_s\}_N$ and
$\{UB_1U^\ast,\dots,UB_tU^\ast\}_N$ are asymptotically free of
second order.
\end{theorem}

\begin{proof}
The proof is a repetition of the proof of Theorem \ref{main:free}
except that we cannot assume that $\{A_1, \dots , A_s, B_1, \dots ,
B_t\}_N$ has a second order limit distribution. Instead we have the
independence of the $A_i$ from the $B_i$'s and a special
$\epsilon$. So we shall only indicate how the proof has to be
modified. 

The first order freeness follows as in the proof of Theorem
\ref{main:free}. In the proofs of Theorems \ref{main},
\ref{main:higherorder}, and \ref{main:free} the cumulants we need
are all sums over $S_{2n}^{(\epsilon)}$ for various $n$'s. Now we
have a special form of $\epsilon$, namely $\epsilon_i =
(-1)^{i+1}$. Thus if $\pi \in S_{2n}^{(\epsilon)}$ then $\pi$ takes
even numbers to odd numbers and vice versa. Since the same applies
to any of the $\gamma$'s, we have that $\gamma \pi^{-1}$ takes even
numbers to even numbers and odd numbers to odd numbers. Thus
the orbits of $\gamma\pi^{-1}$ consist either of all odd numbers
or of all even numbers. Hence if $P_1, \dots , P_{n}$ are words in
$A_1, \dots , A_s$ and $Q_1, \dots , Q_{n}$ are words in $B_1,
\dots , B_t$, then by the independence of the $A$'s and the $B$'s 
\begin{multline*} 
\EE\big( \Tr_{\gamma\pi^{-1}}(P_1, Q_1, \dots , P_{n}, Q_n) \big) \\
= \EE\big( \Tr(W_1) \cdots \Tr(W_j) \big) \cdot \EE\big(
\Tr(W_{j+1}) \cdots \Tr(W_k) \big) 
\end{multline*}
where $k = \#(\gamma\pi^{-1})$ and $W_i$ for $1 \leq i \leq j$ is a
word only in $A$'s and for $j+1 \leq i \leq k$ is a word only
in $B$'s. This means that as far as the asymptotic behaviour of
$\EE\big( \Tr_{\gamma\pi^{-1}}(P_1, Q_1, \dots , P_{n}, Q_n) \big)$
is concerned we may assume that $\{A_1, \dots , A_s, B_1, \dots ,
B_t\}_N$ has a second order limit distribution. Hence our claim
follows from Theorem \ref{main:free} 
\end{proof}

We say that a tuple $\{B_1,\dots,B_s\}$ of $N\times N$-random
matrices is $\UU(N)$--invariant if for every $U\in\UU(N)$ the joint
probability distribution of the random matrices $\{B_1,\dots,B_s\}$
coincides with the joint probability distribution of the random
matrices $\{UB_1U^\ast,\dots,UB_sU^\ast\}$.

\begin{corollary}\label{cor:3.14}
Let $\{A_1,\dots,A_s\}_N$ be a sequence of $N\times N$-random
matrices which has a second order limit distribution and let $\{B_1,
\ab \dots, \ab B_t\}_N$ be a sequence of $\UU(N)$--invariant
$N\times N$-random matrices which has a second order limit
distribution. Furthermore assume that the matrices
$\{A_1,\dots,A_s\}_N$ and the matrices $\{B_1,\dots,B_t\}_N$ are
independent. Then the sequences $\{A_1,\dots,A_s\}_N$ and
$\{B_1,\dots,B_t\}_N$ are asymptotically free of second order.
\end{corollary}

\section{Failure of Universality for multi-matrix models}
In this section we want to make the meaning of second order freeness
for fluctuations of random matrices more explicit and relate this
with the question of universality of such fluctuations. There has
been much interest in global fluctuations of random matrices, in
particular, since it was observed that for large classes of
one-matrix models these fluctuations are universal. In the physical
literature this observation goes at least back to Politzer \cite{P},
culminating in the paper of Ambj\o{}rn et al. \cite{AJM}, whereas a
proof on the mathematical level of rigour is due to Johansson
\cite{J}. Universality for one-matrix models lets one expect
that one would also have such universality for multi-matrix
models. Indeed, this expectation was one of the starting points
of our investigations. However, our machinery around second order
freeness shows that such universality is {\it not}\/ present in
multi-matrix models.

Before we address multi-matrix models let us first recall the
relevant result of Johansson \cite{J}. We consider Hermitian
$N\times N$-random matrices $A=(a_{ij})_{i,j=1}^N$ equipped with the
probability measure
\begin{equation}\label{ensemble}
d\mu_N(A)=\frac 1{Z_N}\exp\bigl\{-N\Tr[P(A)]\bigr\}dA,
\end{equation}
where
$$dA=\prod_{1\leq i<j\leq N}d\,\text{Re}\, a_{ij}\, d\,\text{Im}\, a_{ij}
\prod_{i=1}^N da_{ii}.$$ Here, $P$ is a polynomial in one variable,
which we will address as ``potential" in the following, and $Z_N$ is
a normalization constant. If one restricts to a special class
$\mathcal{V}$ of potentials $P$ (the most important condition being
that the limit eigenvalue distribution has a single interval as
support -- which we normalize in the following to $[-2,2]$) then
Johansson proved the following universality of fluctuations for this
class: Consider the sequence of $N\times N$-random matrices
$\{A_N\}_N$ given by (\ref{ensemble}). Then this sequence has a
second order limit distribution which can be described as follows:
\begin{enumerate}
\item
We have
$$\lim_{N\to\infty}
\kk_1\{\tr[A^n_N]\}=\int t^nd\nu_P(t)\qquad (n\in\NN),$$ where
$\nu_P$ is a probability measure on $\RR$ (``limiting eigenvalue
distribution") which is given as the solution of the singular
integral equation
\begin{equation}\label{integral}
\int\frac{d\nu_P(t)}{t-s}=-\frac 12 P'(s)\qquad \text{for all
$s\in\text{supp } \nu_P$}.
\end{equation}
\item
Let $T_n$ ($n\in\NN$) be the Chebyshev polynomials of first kind
(which are the orthogonal polynomials with respect to the arcsine
law on $[-2,2]$). Then
$$\lim_{N\to\infty}
\kk_2\{\Tr[T_n(A_N)],\Tr[T_m(A_N)]\}=\delta_{mn}\cdot n.$$
\end{enumerate}
Whereas the limiting eigenvalue distribution $\nu_P$ depends on the
form of the potential $P$, the fluctuations are the same for all
potentials in the class $\mathcal{V}$ -- they are always
diagonalized by the same polynomials $\{T_n\mid n\in\NN\}$. Note
that the most prominent example for the considered class of random
matrices is given by $P(A)=A^2$, which corresponds to the Gaussian
random matrix ensemble. So one can phrase this universality also in
the way that all considered random matrices have the same
fluctuations as Gaussian random matrices (but their eigenvalue
distributions are of course different from Wigner's semi-circle
law).

Let us now consider multi-matrix models. For notational convenience
we restrict to the case of two-matrix models. Take now a polynomial
$P(A,B)$ in two non-commuting variables and consider pairs of
hermitian $N\times N$-matrices $A=(a_{ij})_{i,j=1}^N$ and
$B=(b_{ij})_{i,j=1}^N$ equipped with the probability measure
\begin{equation}\label{ensemble-two}
d\mu_N(A,B)=\frac 1{Z_N}\exp\{-N\Tr[P(A,B)]\}dA dB,
\end{equation}
where $Z_N$ is a normalization constant. The simplest ensemble of
this kind is the case of two independent Gaussian random matrices
which corresponds to the choice $P(A,B)=A^2+B^2$. The above
mentioned universality result for the one-matrix case lets one
expect that one might also have universality for multi-matrix
ensembles which are close to the ensemble of independent Gaussian
random matrices. However, we will now show that even restricted to
potentials of the form $P(A,B)= P_1(A)+P_2(B)$ we do not have
universal fluctuations.

Let us first observe that our concept of second order freeness tells
us how to diagonalize fluctuations. We spell this out in the
following theorem which is an easy consequence of Definition
\ref{free} of second order freeness.

\begin{theorem}\label{meaning}
Let $\{A_N\}_N$ and $\{B_N\}_N$ be two sequences of $N\times
N$-random matrices which are asymptotically free of second order.
Let $\{Q_n^A\mid n\in\NN\}$ and $\{Q_n^B\mid n\in\NN\}$ be the
orthogonal polynomials for the limiting eigenvalue distribution of
$A_N$ and $B_N$, respectively, determined by the requirements that
$Q_n^A$ and $Q_n^B$ are polynomials of degree $n$ and that
$$\lim_{N\to\infty}\kk_1\{\tr[Q_n^A(A_N)Q_m^A(A_N)]\}=\delta_{nm},$$
$$\lim_{N\to\infty}\kk_1\{\tr[Q_n^B(B_N)Q_m^B(B_N)]\}=\delta_{nm}.$$
Then the fluctuations of mixed traces in $A_N$ and $B_N$ are
diagonalized by cyclically alternating products of $Q_n^A$ and
$Q_n^B$ and the covariances are given by the number of cyclic
matchings of these products:
\begin{multline*}
\lim_{N\to\infty}
\kk_2\bigl\{\Tr[Q^A_{i(1)}(A_N)Q^B_{j(1)}(B_N)\cdots
Q^A_{i(m)}(A_N)Q^B_{j(m)}(B_N)],\\
\Tr[Q^B_{l(n)}(B_N) Q^A_{k(n)}(A_N) \cdots
Q^B_{l(1)}(B_N) Q^A_{k(1)}(A_N)]\bigr\}\\
=\delta_{mn} \cdot\#\{r\in\{1,\dots,n\}\!\mid\! i(s)=k(s+r),
j(s)=l(s+r)\,\forall s=1,\dots,n\}
\end{multline*}
\end{theorem}

To come back to our problem of universality for multi-matrix models
we only have to observe that Corollary \ref{cor:3.14} tells us that
we have asymptotic freeness of second order if we choose a potential
of the form $P(A,B)=P_1(A)+P_2(B)$.

\begin{theorem}
Let $P(A,B)=P_1(A)+P_2(B)$ where $P_1$ and $P_2$ are polynomials
from the class $\mathcal{V}$. Consider the two-matrix model
$\{A_N, \ab B_N \}_N$ given by the probability measure
(\ref{ensemble-two}). Let $\nu_1$ and $\nu_2$ be the limiting
eigenvalue distribution for $P_1$ and $P_2$, respectively (as
described in Eq. (\ref{integral})) and denote by $\{Q^1_n\mid
n\in\NN\}$ and $\{Q^2_n\mid n\in\NN\}$ the respective orthogonal
polynomials. Then the sequence $\{A_N,B_N\}_N$ has a second order
limit distribution given by $(\cA,\ff_1,\ff_2)$ and $a,b\in\cA$
which can be described as follows:
\begin{enumerate}
\item
$\ff_1$ is the free product of $\nu_1$ and $\nu_2$.
\item
$\ff_2$ is diagonalized by the following collection of polynomials:
$$\{T_n(a)\mid n\in\NN\}$$
$$\{T_n(b)\mid n\in\NN\}$$
$$\{Q^1_{i(1)}(a)Q^2_{j(1)}(b)\cdots
Q^1_{i(n)}(a)Q^2_{j(n)}(b) \},$$ 
where in the last set we choose one representative from each cyclic
equivalence class, i.e. for all $n\in\NN$,
$n$-tuples $\bigl((i(1),j(1)), \ab \dots , \ab
(i(n), \ab j(n))\bigr)$   which are
different modulo cyclic rotation
\end{enumerate}
\end{theorem}

\begin{proof}
Note that the additive form of the potential $P(A,B)=P_1(A)+P_2(B)$
means that $A_N$ and $B_N$ are independent, $A_N$ is a one-matrix
ensemble corresponding, via (\ref{ensemble}), to a potential $P_1$,
and $B_N$ is a one-matrix ensemble corresponding to a potential
$P_2$. Thus the statement about the diagonalizing polynomials in
either only $A_N$ or in only $B_N$ is just Johansson's result. For
getting the statement about the diagonalizing polynomials in both
$A_N$ and $B_N$ we have to note that the random matrices $A_N$ (and
also $B_N$) are $\UU(N)$-invariant, thus Corollary \ref{cor:3.14}
implies that $\{A_N\}_N$ and $\{B_N\}_N$ are asymptotically free of
second order. Hence we can apply Theorem \ref{meaning} 
above.\end{proof}

Note that whereas the polynomials in only one of the matrices are
universal (namely equal to the Chebyshev polynomials $\{T_n\}_n$), the
polynomials which involve both matrices are not universal but depend
on the eigenvalue distributions $\nu_1$ and $\nu_2$. Since the
latter vary with the potentials $P_1$ and $P_2$, the polynomials
$Q^1_n$ and $Q^2_n$, and thus also the alternating products in them,
depend on the choice of $P$. Thus we can conclude that even within
the very restricted class of $P$'s of the form
$P(A,B)=P_1(A)+P_2(B)$ we have no universality of global
fluctuations in multi-matrix models.

\end{document}